\documentclass[12pt]{article}
\usepackage{amsfonts,amssymb}

\def\C{\mathbb{C}}
\def\N{\mathbb{N}}
\def\Z{\mathbb{Z}}

\def\bq{ \begin{equation} }
\def\eq{ \end{equation} }
\def\ben{ \begin{eqnarray} }
\def\en{ \end{eqnarray} }

\def\frac#1#2{{#1\over #2}}

\def\on#1#2{\mathop{\vbox{\ialign{##\crcr\noalign{\kern2pt}
$\scriptstyle{#2}$\crcr\noalign{\kern2pt\nointerlineskip}
\kern-2pt$\hfil\displaystyle{#1}\hfil$\crcr}}}\limits}

\begin{document}

\baselineskip=15pt
\vspace{1cm} \noindent {\LARGE \textbf {Algebraic structures
connected with pairs of compatible associative algebras}} \vskip1cm
\hfill
\begin{minipage}{13.5cm}
\baselineskip=15pt {\bf
 Alexander Odesskii ${}^{1,} {}^{2}$ and
 Vladimir Sokolov ${}^{1}$} \\ [2ex]
{\footnotesize ${}^1$ Landau Institute for Theoretical Physics,
Moscow, Russia
\\
${}^{2}$  School of Mathematics, The University of
Manchester, UK}\\
\vskip1cm{\bf Abstract}

We study associative multiplications in semi-simple associative
algebras over $\C$ compatible with the usual one or, in other words,
linear deformations of semi-simple associative algebras over $\C$.
It turns out that these deformations are in one-to-one
correspondence with representations of certain algebraic structures,
which we call $M$-structures in the matrix case and $PM$-structures
in the case of direct sums of several matrix algebras. We also
investigate various properties of $PM$-structures, provide numerous
examples and describe an important class of $PM$-structures. The
classification of these $PM$-structures naturally leads to affine
Dynkin diagrams of $A$, $D$, $E$-type.

\end{minipage}

\vglue1cm \textbf{Address}: Landau Institute for Theoretical
Physics, Kosygina 2, 119334, Moscow, Russia

\textbf{E-mail}: odesskii@itp.ac.ru, \, sokolov@itp.ac.ru \newpage

\centerline{\Large\bf Introduction}
\medskip
Two associative algebras with multiplications $\star$ and $\circ$
 defined on
the same finite dimensional vector space $\bf V$ are said to be {\it
compatible} if the multiplication
\begin{equation} \label{pensil}
a \bullet b= a\star b+ \lambda \,a\circ b
\end{equation}
is associative for any constant $\lambda$. The multiplication
$\bullet$ can be regarded as a deformation of the multiplication
$\star$ linear in parameter $\lambda$.

The description of pairs of compatible associative products seems to
be an interesting mathematical problem on its own. Moreover, the
approach to integrable systems based on the concept of compatible
Poisson structures via Lenard-Magri scheme \cite{magri} provides
further  motivation for investigation of compatible associative
multiplications.

Recall that two Poisson brackets are said to be {\it compatible} if
any linear combination of these brackets is a Poisson bracket. It is
well-known that the formula $ \{x_i,\, x_j \}=c^k_{ij} x_k, \quad
i,j=1,\dots,N $ defines a linear Poisson structure iff $c^k_{ij}$
are structural constants of a Lie algebra. The compatibility of two
such structures is equivalent to the compatibility of the
corresponding Lie brackets.  Various applications of compatible Lie
brackets in the integrability theory can be found in
\cite{golsok1,bolsbor, golsok2,golsok3,reysem}.

Suppose now that we have two compatible associative algebras with
multiplications $\star$ and $\circ$ defined on the same finite
dimensional vector space $\bf V$. We can construct immediately two
compatible Lie brackets by the usual formulas $[a,b]_1=a\star
b-b\star a$ and $[a,b]_2=a\circ b-b\circ a$ and hence, two
compatible linear Poisson structures.

Moreover, for any $n\in \N$ we can construct two compatible
associative algebras in the space $Mat_n(\bf V),$ which is the space
of $n\times n$ matrices with entries from $\bf V$. Therefore, for
each $n$ we have a pair of compatible Poisson structures in the
linear space of dimension $n^2\dim \bf V$. Note  that even if both
associative algebras on $\bf V$ are commutative we have nontrivial
Poisson structures on the space $Mat_n(\bf V)$ for $n>1$. In terms
of coordinates, if $\{e_i, i=1,...N\}$ is a basis of $\bf V$ and
$e_i\star e_j=p_{i,j}^ke_k, \quad e_i\circ e_j=q_{i,j}^k e_k,$ then
for each $n$ we have two compatible Poisson structures given, in
coordinates $\{f_{i,l,m}, \quad i=1,...,N, \quad l,m=1,...n\}$, by
the formulas
$$
\{f_{i,l_1,m_1},f_{j,l_2,m_2}\}_1=\delta_{m_1,l_2}p_{i,j}^kf_{k,l_1,m_2}-
\delta_{m_2,l_1}p_{j,i}^kf_{k,l_2,m_1}
$$
and
$$\{f_{i,l_1,m_1},f_{j,l_2,m_2}\}_2=\delta_{m_1,l_2}q_{i,j}^kf_{k,l_1,m_2}-
\delta_{m_2,l_1}q_{j,i}^kf_{k,l_2,m_1}.
$$
Note that these Poisson structures are invariant with respect to the
action of the group $GL_n(\C)$ on the space $Mat_n(\bf V)$ by
conjugations. Therefore, for any two functions invariant with
respect to this action their Poisson bracket is also invariant.
Since any invariant function can be written in terms of traces of
matrix polynomials, we see that a bracket of two traces can be also
written in terms of traces. This leads us to bi-hamiltonian
structures for the so-called nonabelian integrable systems in the
sense of \cite{sokmih}.

Another motivation for the study of compatible associative algebras
can be found in \cite{marmo}.

In this paper we assume that the associative algebra over field $\C$
with multiplication $\star$ is semi-simple. In other words, this
algebra is a direct sum of matrix algebras over $\C$ \cite{alg}. It
turns out that in this case  multiplications $\circ$ compatible with
$\star$ are in one-to-one correspondence to representations of
special infinite-dimensional associative algebras. The simplest
finite-dimensional version of such an algebra can be described as
follows. Let ${\cal A}$ and ${\cal B}$ be associative algebras of
the same dimension $p$ with bases $A_{1},\dots, A_{p}$ and
$B^{1},\dots, B^{p}$ and structural constants $\phi^i_{j,k}$ and
$\psi^{\alpha,\beta}_{\gamma},$ correspondingly. Suppose that the
structural constants satisfy the following identities:
$$
 \phi^s_{j,k}\psi_s^{l,i}=\phi^l_{s,k}\psi_j^{s,i}+\phi^i_{j,s}\psi_k^{l,s}, \qquad 1\le i,j,k,l \le p.
$$
Then the algebra  of dimension $2p+p^{2}$ with the basis $A_{i},
B^{j},A_{i} B^{j}$ and relations
$$
B^i A_j=\psi_j^{k,i} A_k+\phi^i_{j,k} B^k
$$
is associative.

An invariant description of such a construction can be given as
follows. Suppose that we have two associative algebras $\cal A$ and
$\cal B$, a non-degenerate pairing $\cal A\times\cal B\to\C$ and
structures of right $\cal A$-module and left $\cal B$-module on the
space $\cal A\oplus\cal B$ commuting with each other. Assume also
that $\cal A$ acts in this module by right multiplication on itself
and $\cal B$ acts by left multiplication on itself. Extend our
pairing to the space $\cal A\oplus\cal B$ by the formulas
$(a_1,a_2)=(b_1,b_2)=0$ for $a_1,a_2\in\cal A$, $b_1,b_2\in\cal B$
and assume that it is invariant under the action of $\cal A$ and
$\cal B$: $(va_1,a_2)=(v,a_1a_2)$ and $(b_1,b_2v)=(b_1b_2,v)$ for
 $v\in\cal A\oplus\cal B$.
In this situation one can define a natural structure of an
associative algebra on the space $\cal A\oplus\cal B\oplus (\cal
A\otimes\cal B)$ compatible with our module structures. This means
that the action of $\cal A$ on the algebra by right multiplication
restricted to $\cal A\oplus\cal B$ coincides with the module action
and the same property is valid for the action of $\cal B$ by left
multiplication.

Algebras considered in this paper are more complicated. Namely,
algebras $\cal A$ and $\cal B$ have common unity. Instead of their
direct sum, we construct a linear space of the same dimension but
with one-dimensional ``defect'': $\cal A$ and $\cal B$ are
intersected by the linear span of the unity but we add one more
element, which in some sense dual to the unity. We assume the
existence of a non-degenerate pairing and structures of a right
$\cal A$-module and a left $\cal B$-module on this space with
properties similar to described above. A linear space with these
structures and the corresponding associative algebra are called {\it
$M$-structure} and {\it $M$-algebra}. In turns out that $M$-algebra
is infinite-dimensional over $\C$ but finite dimensional over the
subalgebra generated by a special central element $K$.

The main result of this paper is the following: there is a
one-to-one correspondence between $n$-dimensional representations
(that should be non-degenerate in some sense) of $M$-algebras
 and associative products in
$Mat_n$ compatible with the usual matrix product. In other words, to
find all associative products in matrix algebras compatible with the
usual one, we should describe $M$-structures and for each
$M$-structure classify finite-dimensional representations of the
corresponding $M$-algebra.

To describe the compatible products for the algebra
$Mat_{n_1}\oplus...\oplus Mat_{n_m}$ we introduce {\it
$PM$-algebras}, which are  generalizations of $M$-algebras. Roughly
speaking, a $PM$-algebra looks like the algebra of $m\times m$
matrices with entries being elements of some $M$-algebra.

This paper is organized as follows. In  Section {\bf 1}, we collect
some general facts about compatible associative multiplications. The
first result of that Section is standard and based on the
deformation theory of associative algebras. Namely, we show that if
the algebra with respect to multiplication $\star$ is rigid (which
holds in semi-simple case), then there exists a linear operator
$R:\bf V\to\bf V$ such that the  multiplication $\circ $ is of the
form
\begin{equation} \label{mult2}
X \circ Y =R(X) \star Y+X \star R(Y)-R(X \star Y).
\end{equation}
We also provide several examples of compatible multiplications. At
the end of Section {\bf 1} we give a construction of $m+1$ pairwise
compatible associative multiplications on the space ${\bf V}\otimes
F_m$ provided that we have two compatible associative
multiplications on the space $\bf V$. Here $F_m$ is the space of
polynomials in one variable of degree less then $m$.

In Section {\bf 2}, we consider multiplications compatible with the
standard matrix product in  $Mat_n$. In Subsection {\bf 2.1} we
study admissible operators $R$ written in the form
\begin{equation}\label{Rmat}
R(x)=a_1 \,x \,b^1+...+a_p \,x\, b^p+c\, x
\end{equation}
with $p$ being smallest possible. It turns out that
$a_1,...,a_p,b^1,...,b^p,c$ should be generators of a representation
of an $M$-structure. In Subsection {\bf 2.2}, we propose an
invariant definition of $M$-structures and $M$-algebras and study
their properties.  In
 Subsection {\bf 2.3} we describe $M$-algebras
in the special case when the algebra $\cal A$ is commutative
semi-simple (that is, isomorphic to $\C\oplus...\oplus\C$).

Section {\bf 3} is devoted to a generalization of the results from
the previous section to the case of the algebra
$Mat_{n_1}\oplus...\oplus Mat_{n_m}$. All results and proofs are
similar to the ones from Section {\bf 2}. In Subsection {\bf 3.1} we
study possible operators $R$, and in Subsection {\bf 3.2} give an
invariant definition of the corresponding algebraic structures.

In  Section {\bf 4} we describe all $PM$-structures with semi-simple
algebras $\cal A$ and $\cal B$. It turns out that such
$PM$-structures are related to Cartan matrices of affine Dynkin
diagrams of the $\tilde A_{2 k-1},$ $\tilde D_{k},$ $\tilde E_{6},$
$\tilde E_{7},$ and $\tilde E_{8}$-type.

In  Conclusion we discuss some implications of our results and
possible directions of further research.

\section{Compatible associative multiplications}

Suppose that we have an associative multiplication $\star$ defined
on a finite dimensional vector space $\bf V$ such that $\bf V$ is a
semi-simple algebra with respect to this multiplication. The
following classification problem arises: to describe all possible
associative multiplications $\circ$ on a vector space ${\bf V},$
compatible with a given semi-simple multiplication $\star$. Since
any semi-simple associative algebra is rigid, the multiplication
(\ref{pensil}) is isomorphic to $\star$ for almost all values of the
parameter $\lambda$. This means that there exists a formal series of
the form
\begin{equation}
A_{\lambda}=1+ R \ \lambda+ S \ \lambda^2 + \cdots,   \label{RS}
\end{equation}
where the coefficients are linear operators on {\bf V}, such that
\begin{equation}
A_{\lambda}^{-1}\Big(A_{\lambda}(X) \star A_{\lambda}(Y)\Big)=  X
\star Y+\lambda \ X \circ Y. \label{sog2}
\end{equation}
Equating the coefficients of $\lambda$ in (\ref{sog2}), we obtain
the formula (\ref{mult2}). It is easy to see that the transformation
\begin{equation}
     R\longrightarrow R + ad_{\star} a,
\label{ad}
\end{equation}
does not change the multiplication $\circ$ for any $a\in {\bf V},$
where $ad_{\star}a$ is a linear operator $v\to a\star v-v\star a$.

Comparing the coefficients of $\lambda^2$ in (\ref{sog2}), we get
the following identity
\begin{equation}
\begin{array}{l} \label{yangRS}
R\big(R(X)\star Y + X\star R(Y)\big) - R(X)\star R(Y) - R^2(X \star
Y)
\\[3mm]
=S(X) \star Y+X \star S(Y)-S(X \star Y),
\end{array}
\end{equation}
for any $X,Y\in {\bf V}$. It is not difficult to prove that if
$(\ref{yangRS})$ holds for some operators $R$ and $S$ then the
multiplication  (\ref{mult2}) is associative and compatible with
$\star$. Under transformation (\ref{ad}) the operator $S$ is
changing as follows
$$
 S\longrightarrow S + ad_{\star} a \circ R+\frac{1}{2} (
 ad_{\star}
a)^2.
$$

In the important special case $S=0,$ we have
\begin{equation}
 \label{yangRR}
R\,\Big(R(X)\star Y + X\star R(Y)\Big) - R(X)\star R(Y) - R^2(X
\star Y) =0.
\end{equation}
In the paper \cite{marmo} some examples of such $R$-operators have
been found.

{\bf Definition.} We call operators $R$ and $R^{\prime}$ equivalent
if $R-R^{\prime}=ad_{\star}a$ for some $a\in \bf V$.

\noindent It is known that any derivative of semi-simple algebra has
the form $ad_{\star}a$ for some $a\in \bf V$. Therefore, the formula
(\ref{mult2}) gives the same multiplications for operators $R$ and
$R^{\prime}$ if and only if these operators are equivalent.

{\bf Example 1.1.} Let $a$ be an arbitrary element of $\bf V$ and
$R$ be the operator of left multiplication by $a$ with respect to
$\star$. Then $R$ satisfies (\ref{yangRR}) and the corresponding
multiplication $X\circ Y=X\star a\star Y$ is associative and
compatible with $\star$.

{\bf Example 1.2.} Suppose that $\star$ is the standard matrix
product in ${\bf V}=Mat_2\,$, $a,b\in {\bf V}$, then the product
$$
X\circ Y=(a X-X a)\,(b Y-Y b)
$$
is associative and compatible with the standard one. The
corresponding operator $R$ is given by $R(X)=a\, (X b-b X).$ If
$\mbox{Det}\,a =0 $, then the operator $R$  satisfies
(\ref{yangRR}). The Example 1 from the paper \cite{marmo}
corresponds to the special case of the Example 1.2 where the
matrices $a$ and $b$ are diagonal.

The following statement can be verified straightforwardly.

{\bf Proposition 1.1.}  The Examples 1.1 and 1.2 describe all
associative multiplications compatible with the matrix product in
$Mat_2.$

{\bf Example 1.3.} Let ${\bf e}_{1},\dots,{\bf e}_{m}$ be a basis in
$\bf V$ and the multiplication $\star$ is given by
\begin{equation}\label{alg}
{\bf e}_{i}\star {\bf e}_{j}=\delta^{i}_{j} {\bf e}_{i}.
\end{equation}
The algebra thus defined is commutative and semi-simple.  Suppose
the entries $r_{ij}$ of the matrix $R$ satisfy the following
relations:
$$
\sum_{k=1}^{m} r_{ki}=q_{0}, \qquad \hbox{and} \quad \qquad r_{ik}
r_{jk}=r_{ij} r_{jk}+r_{ji} r_{ik} \quad \hbox{for} \quad i\ne j \ne
k\ne i,
$$
where $q_{0}$ is an arbitrary constant. The generic solution of this
system of algebraic equations is given by
$$
r_{ii}=q_{0}- \sum_{k\ne i} r_{ki}, \qquad    r_{ij}=\frac{q_{i}
p_{j}}{p_{i}-p_{j}}, \qquad i\ne j,
$$
where $p_{i},  q_{j}$ are arbitrary constants. The formula
(\ref{mult2}) defines a multiplication
$$
{\bf e}_{i}\circ {\bf e}_{j}=r_{ij} {\bf e}_{j}+r_{ji} {\bf
e}_{i}-\delta^{i}_{j}\sum_{k=1}^{m} r_{ik} {\bf e}_{k}
$$
compatible with $\star$. Since this multiplication is linear with
respect to the parameters $q_{i},$ we have got a family of $m+1$
pairwise compatible associative multiplications. This family can be
described in a different way in terms of the generating function
$$
 {\bf E}(u)={\bf e}_1+u\, {\bf e}_2+\dots+u^{m-1}\, {\bf e}_{m}.
$$
Let $ q(u)=a_0+u\,a_1+\dots+u^m\, a_{m}$ be an arbitrary polynomial
of degree $n$. Define a multiplication of the generating functions
by the formula
\begin{equation}\label{mults}
{\bf E}(u)\,{\bf E}(v)=\frac{u q(v)}{u-v}\,{\bf E}(u)+\frac{v
q(u)}{v-u} \,{\bf E}(v).
\end{equation}
It is easy to verify that (\ref{mults}) yields an associative
multiplication between $ {\bf e}_1,...,{\bf e}_{m}$ linear with
respect to the parameters $a_0,\dots,a_m.$ Let $b_1,...,b_m$ be
roots of $q(u)$ and assume that these roots are pairwise distinct.
Then $\tilde{\bf e}_{i}=b_{i} q'(b_{i})\,{\bf E}(b_i)$ form a basis,
in which this multiplication is given by (\ref{alg}).

The formula (\ref{mults}) admits the following generalization. Let
$\bf V$ be a finite dimensional vector space with two compatible
associative multiplications $\star$ and $\circ$. Let $F_m$ be a
vector space of polynomials in one variable $t$ with degree less
then $m$. We are going to construct $m+1$ pairwise compatible
associative multiplications on the space ${\bf V}\otimes F_m$. For a
vector $x\in\bf V$ we denote by $x_i$ the element $x\otimes
t^i\in{\bf V}\otimes F_m$. Denote by $x(u)$ the following polynomial
in $u$ with values in the space ${\bf V}\otimes F_m$:
$$x(u)=x_0+x_1u+...+x_{m-1}u^{m-1}.$$ Let us also fix an arbitrary
polynomial $q(u)\in\C[u]$ of degree $m$.

{\bf Theorem 1.1.} The formula
\begin{equation}\label{mults1}
x(u)y(v)=\frac{q(u)}{u-v}((x\star y)(v)+v(x\circ
y)(v))+\frac{q(v)}{v-u}((x\star y)(u)+u(x\circ y)(u))
\end{equation}
defines an associative multiplication on the linear space ${\bf
V}\otimes F_m$. Here $x,y\in\bf V$ are arbitrary vectors.

{\bf Proof.} It is clear that both r.h.s. and l.h.s. of
(\ref{mults1}) are polynomials in $u,v$ of degree $m-1$ with values
in ${\bf V}\otimes F_m$. Therefore, the formula defines a product in
this space. Associativity of this product can be easily checked by
direct calculation.

Note that the formula (\ref{mults1}) defines the product which
linearly depends on the polynomial $q(u)$ of degree $m$. Therefore,
we have $m+1$ pairwise compatible associative multiplications on the
space ${\bf V}\otimes F_m$.

{\bf Remark 1.} If $\bf V=\C$ with trivial pair $1\star 1=0, 1\circ
1=1$, then this construction gives the Example 1.3 (see
(\ref{mults})).

{\bf Remark 2.} Let $b_1,...,b_m$ be roots of $q(u)$ and assume that
these roots are pairwise distinct. One can check that the algebra
${\bf V}\otimes F_m$ with respect to the multiplication
(\ref{mults1}) is isomorphic to a direct sum of $m$ components.
Moreover, the $i$th component is isomorphic to $\bf V$ with respect
to the product $x\bullet y=x\star y+b_ix\circ y$. This is a direct
consequence of the formula (\ref{mults1}). In particular, if $\bf V$
is semi-simple for generic linear combination of $\star$ and
$\circ$, and the roots $b_1,...,b_m$ are also generic, then ${\bf
V}\otimes F_m$ is isomorphic to direct sum of $m$ copies of $\bf V$.

\section{Matrix case}

\subsection{Construction of the second product}

Consider now the case where the algebra is isomorphic to $Mat_n$
with respect to the first product. Any linear operator $R$ on the
space $Mat_n$ may be written in the form $R(x)=a_1 x b^1+...+a_l x
b^l$ for some matrices $a_1,...,a_l,b^1,...,b^l$. Indeed, the
operators $x\to e_{i,j} x e_{i_1,j_1}$ form a basis in the space of
linear operators on $Mat_n$.

It is convenient to represent the operator $R$ from the formula
(\ref{mult2}) in the form (\ref{Rmat})
 with $p$  being smallest possible in the class of
equivalence of $R$. This means that the matrices $\{a_1,...,a_p,1\}$
are linear independent as well as the matrices $\{b^1,...,b^p,1\}$.
According to (\ref{mult2}), the second product has the following
form
\begin{equation}\label{sekprod}
x\circ y=a_i\, x\, b^i\, y+x \,a_i\, y\, b^i-a_i\, x y\, b^i+x\,
c\, y.
\end{equation}

Note that we have the following transformations, which do not change
the class of equivalence of $R$. The first family of such
transformations is
$$a_i\to a_i+u_i, \qquad b^i\to b^i+v^i, \qquad c\to
c-u_ib^i-v^ia_i-u_iv^i
$$
 for any constants $u_1,...,u_p,v^1,...,v^p$ and
the second one is
$$
a_i\to g_i^ka_k, \qquad b^i\to h_k^ib^k, \qquad c\to c,
$$
 where
$g_i^kh_k^j=\delta_i^j$. This means that we can regard $a_i$ and
$b^i$ as bases in dual vector spaces.

{\bf Theorem 2.1.} The multiplication $\circ$ given by the formula
(\ref{sekprod}) is an associative product on the space $Mat_n$ if
and only if there exist tensors $\phi^k_{i,j}, \mu_{i,j},
\psi^{i,j}_k, \lambda^{i,j}, t^i_j$ such that the following
relations hold:
\begin{equation}\label{(!)}
a_ia_j=\phi_{i,j}^ka_k+\mu_{i,j},\qquad
b^ib^j=\psi_k^{i,j}b^k+\lambda^{i,j},
\end{equation}
\begin{equation}\label{(!!)}
b^ia_j=\psi_j^{k,i}a_k+\phi^i_{j,k}b^k+t^i_j+\delta^i_jc,
\end{equation}
\begin{equation}\label{(!!!)}
b^ic=\lambda^{k,i}a_k-t^i_kb^k-\phi^i_{k,l}\psi^{l,k}_sb^s-\phi^i_{k,l}\lambda^{l,k},\quad
ca_j=\mu_{j,k}b^k-t^k_ja_k-\phi^s_{k,l}\psi^{l,k}_ja_s-\mu_{k,l}\psi^{l,k}_j
\end{equation}
Moreover, the tensors $\phi^k_{i,j}, \mu_{i,j}, \psi^{i,j}_k,
\lambda^{i,j}, t^i_j$ satisfy the properties
\begin{equation}\label{(!!!!)}
\phi^s_{j,k}\phi^i_{s,l}+\mu_{j,k}\delta^i_l=\phi^i_{j,s}\phi^s_{k,l}+\delta^i_j\mu_{k,l},\qquad
 \phi^s_{j,k}\mu_{i,s}=\phi^s_{i,j}\mu_{s,k},
 \end{equation}
 \begin{equation}\label{(!!!!!)}
 \psi^{i,j}_s\psi^{s,k}_l+\delta^k_l\lambda^{i,j}=\psi^{j,k}_s\psi^{i,s}_l+\delta^i_l\lambda^{j,k},
 \qquad
 \psi_s^{i,j}\lambda^{s,k}=\psi^{j,k}_s\lambda^{i,s},
 \end{equation}
 $$\phi^s_{j,k}\psi_s^{l,i}=\phi^l_{s,k}\psi_j^{s,i}+\phi^i_{j,s}\psi_k^{l,s}+\delta^l_kt^i_j-
 \delta^i_jt_k^l-\delta^i_j\phi^l_{s,r}\psi^{r,s}_k,$$
 \begin{equation}\label{(!!!!!!)}
 \phi^s_{j,k}t^i_s=\psi_j^{s,i}\mu_{s,k}+\phi^i_{j,s}t_k^s-\delta^i_j\psi_k^{s,r}\mu_{r,s},
 \quad
 \psi^{k,i}_st_j^s=\phi^i_{j,s}\lambda^{k,s}+\psi_j^{s,i}t_s^k-\delta^i_j\phi^k_{s,r}\lambda^{r,s}
 \end{equation}

{\bf Proof.} Associativity $(x\circ y)\circ z=x\circ (y\circ z)$ is
equivalent to the following identity
\begin{equation}
\begin{array}{l}
a_ia_jx(b^jyb^i-yb^jb^i)z+a_ia_jx(yb^j-b^jy)zb^i+x(a_ia_jy-a_iya_j)zb^jb^i+\\[3mm]
a_ix(ya_j-a_jy)zb^jb^i+a_ix(b^ia_jy-yb^ia_j)zb^j+a_ix(a_jyb^jb^i-b^ia_jyb^j+cyb^i-yb^ic)z+\\[3mm]
x(a_jyb^ja_i-a_ia_jyb^j-a_iyc+ca_iy)
zb^i+x(a_iyb^ic-ca_iyb^i)z+a_ix(yc-cy)zb^i=0
\end{array} \label{(propas)}
\end{equation}
From this identity one can readily obtain (\ref{(!)}),
(\ref{(!!)}), (\ref{(!!!)}) using the following

{\bf Lemma 2.1.} Let $p_1xq_1+...+p_lxq_l=0$ for all $x\in Mat_n$.
If $p_1,...,p_l$ are linear independent matrices, then
$q_1=...=q_l=0$. Similarly, if $q_1,...,q_l$ are linear independent
matrices, then $p_1=...=p_l=0$.

Indeed, suppose that some product $a_{i_0}a_{j_0}$ is linearly
independent of $1,a_1,...,a_p$. Since $1,a_1,...,a_p$ are linear
independent by assumption, there exists such a basis in the linear
space spanned by $\{1,a_i,a_ia_j;1\le i,j\le p\}$ that is a subset
of this set and contains the subset
$\{1,a_1,...,a_p,a_{i_0}a_{j_0}\}$. In this basis the coefficient of
$a_{i_0}a_{j_0}$ has the form
\begin{equation} \label{koef}
q_{i,j} \Big((b^j y b^i-y b^jb^i) z+(y b^j-b^j y)z b^i\Big),
\end{equation}
where $q_{i,j}$ are some constants not all equal to zero. Given $y,
z,$ consider the left hand side of (\ref{(propas)}) as a linear
operator applying to the argument $x.$ It follows from Lemma 2.1
that the coefficient (\ref{koef}) is equal to zero. Applying again
Lemma 2.1 to the operator (\ref{koef}) and using the linear
independence of $1,b^1,...,b^p,$ we obtain $q_{i,j}=0$ for all $i$
and $j,$ which is a contradiction. Therefore, all $a_i a_j$ are
linear combinations of $1,a_1,...,a_p$ and, similarly, all $b^i b^j$
are linear combinations of $1,b^1,...,b^p$. This proves (\ref{(!)}).
Substitute these expressions for $a_ia_j$ and $b^i b^j$ to
(\ref{(propas)}) and apply Lemma 2.1 twice. Firstly, we consider the
left hand side of (\ref{(propas)}) as a linear operator with
argument $x$ and take $1,a_1,...,a_p$ for $p_1,...,p_l.$ After that
we regard the same expression as a linear operator with argument $z$
and take $1,b^1,...,b^p$ for $q_1,...,q_l$. As the result, we obtain
the equation
$[y,\,b^ia_j-\psi_j^{k,i}a_k-\phi^i_{j,k}b^k-\delta^i_jc]=0$
equivalent to (\ref{(!!)}) and the following relations
$$\phi^k_{i,j}(b^jyb^i-yb^jb^i)+\lambda^{j,k}(ya_j-a_jy)+a_jyb^jb^i-b^ka_jyb^j+cyb^k-yb^kc=0,$$
$$\psi_k^{j,i}(a_ia_jy-a_iya_j)+\mu_{k,j}(yb^j-b^jy)+a_jyb^ja_k-a_ka_jyb^j-a_kyc+ca_ky=0.$$
Substituting the expressions (\ref{(!)}) and (\ref{(!!)}) for $a_i
a_j$, $b^i b^j$ and $b^j a_i$ into these relations, we get
(\ref{(!!!)}). It can be checked that all these steps are invertible
and (\ref{(propas)}) follows from (\ref{(!)})-(\ref{(!!!)}).

The associativity of the matrix product $a_i a_j a_k$ and the
linear independence of $a_1,...,a_p, 1$ imply (\ref{(!!!!)}).
Similarly, (\ref{(!!!!!)}) follows from the associativity of the
product $b^i b^j b^k$ and the linear independence of $b^1,...,b^p,
1$. The remaining identities are consequences of the associativity
for $b^i a_j a_k$ and $b^i b^j a_k.$

{\bf Remark.} Under conditions (\ref{(!)})-(\ref{(!!!)}), the
operator (\ref{Rmat}) satisfies (\ref{yangRS}) with
$$
S(x)=\mu_{ji} \Big(b^{i} x b^{j}-\psi^{i,j}_k x b^{k} -\lambda^{ij}
x \Big).
$$
In particular, the operator (\ref{Rmat}) satisfies (\ref{yangRR})
iff $\mu_{ji}=0$.

\subsection{$M$-structures and corresponding associative algebras}

In this subsection we describe the algebraic structure underlying
Theorem 2.1.

{\bf Definition.} By weak $M$-structure on a linear space $\cal L$
we mean the following data:
\begin{itemize}
\item  Two subspaces $\cal A$ and $\cal B$ and distinguished
element $1\in\cal A\cap\cal B\subset \cal L$.

\item  A non-degenerate symmetric scalar product $(\cdot, \cdot)$ on the
space $\cal L$.

\item  Associative products $\cal A\times\cal A\to\cal A$ and
$\cal B\times\cal B\to\cal B$ with unity $1$.

\item  A left action $\cal B\times\cal L\to\cal L$ of the algebra
$\cal B$ and a right action $\cal L\times\cal A\to\cal L$ of the
algebra $\cal A$ on the space $\cal L,$ which  commute to each
other.
\end{itemize}

These data should satisfy the following properties:

{\bf 1.} $\dim{\cal A\cap\cal B=\dim\cal L/(\cal A+\cal B)}= 1$. The
intersection of $\cal A$ and $\cal B$ is the one dimensional space
spanned by the unity $1$.

{\bf 2.} The restriction of the action $\cal B\times\cal L\to\cal L$
to the subspace $\cal B\subset \cal L$ is the product in $\cal B$.
The restriction of the action $\cal L\times\cal A\to\cal L$ to the
subspace $\cal A\subset \cal L$ is the product in $\cal A$.

{\bf 3.} $(a_1,a_2)=(b_1,b_2)=0$ and $(b_1b_2,v)=(b_1,b_2v)$,
$(v,a_1a_2)=(va_1,a_2)$ for any $a_1,a_2\in\cal A$, $b_1,b_2\in\cal
B$ and $v\in\cal L$.

It follows from these properties that $(\cdot,\cdot)$ gives a non-
degenerate pairing between ${\cal A}/\C1$ and ${\cal B}/\C1$, so
$\dim\cal A=\dim\cal B$ and $\dim{\cal L}=2\dim{\cal A}$.

For given weak $M$-structure $\cal L$ we can define an algebra
generated by $\cal L$ with natural compatibility and universality
conditions.

{\bf Definition.} By weak $M$-algebra associated with a weak
$M$-structure $\cal L$ we mean an associative algebra $U(\cal L)$
with the following properties:

{\bf 1.} $\cal L\subset U(\cal L)$ and the actions $\cal B\times\cal
L\to\cal L$, $\cal L\times\cal A\to\cal L$ are restrictions of the
product in $U(\cal L)$.

{\bf 2.} For any algebra $X$ with the property {\bf 1} there exist
and unique a homomorphism of algebras $X\to U(\cal L),$ which is
identity on $\cal L$.

It is easy to see that if $U(\cal L)$ exist, then it is unique for
given $\cal L$. Let us describe the structure of $U(\cal L)$
explicitly. Let $\{1,A_1,...,A_p\}$ be a basis of $\cal A$ and
$\{1,B^1,...,B^p\}$ be a dual basis of $\cal B$ (which means that
$(A_i,B^j)=\delta_i^j$). Let $C\in\cal L$ does not belong to the sum
of $\cal A$ and $\cal B$. Since $(\cdot , \cdot)$ is non-
degenerate, we have $(1,C)\ne 0$. Multiplying $C$ by constant, we
can assume that $(1,C)=1$. Adding linear combination of
$1,A_1,...,A_p,B^1,...,B^p$ to $C,$ we can assume that
$(C,C)=(C,A_i)=(C,B^j)=0.$ Such element $C$ is uniquely determined
by choosing basis in $\cal A$ and $\cal B$.

{\bf Lemma 2.2.} The algebra $U(\cal L)$ is defined by the
following relations
\begin{equation}\label{(walg)}
A_iA_j=\phi_{i,j}^kA_k+\mu_{i,j}, \qquad
B^iB^j=\psi_k^{i,j}B^k+\lambda^{i,j}
\end{equation}
\begin{equation}\label{(walg1)}
B^iA_j=\psi_j^{k,i}A_k+\phi^i_{j,k}B^k+t^i_j+\delta^i_jC,
\end{equation}
\begin{equation}\label{(walg2)}
B^iC=\lambda^{k,i}A_k+u_k^iB^k+p^i, \qquad
CA_j=\mu_{j,k}B^k+u^k_jA_k+q_i
\end{equation}
for certain tensors
$\phi_{i,j}^k,\psi_k^{i,j},\mu_{i,j},\lambda^{i,j},u_k^i,p^i,q_i$.

{\bf Proof.} Relations (\ref{(walg)}) just mean that $\cal A$ and
$\cal B$ are associative algebras. Since $\cal L$ is a left $\cal
B$-module and a right $\cal A$-module, the products
$B^iA_j,CA_j,B^iC$ should be linear combinations of the basis
elements $1,A_1,...,A_p,B^1,...,B^p,C$. Applying property {\bf 3} of
weak $M$-structure, we obtain required form of these products. The
universality condition of $U(\cal L)$ shows that this algebra is
defined by (\ref{(walg)}) - (\ref{(walg2)}).

Let us define an element $K\in U(\cal L)$ by the formula $K=A_i
B^i+C$. It is clear that $K$ thus defined does not depend on the
choice of the basis in $\cal A$ and $\cal B$ provided
$(A_i,B^j)=\delta_i^j$, $(1,C)=1$ and $(C,C)=(C,A_i)=(C,B^j)=0$.
Indeed, the coefficients of $K$ are just entries of the tensor
inverse to the form $(\cdot, \cdot)$.

{\bf Definition.} Weak $M$-structure $\cal L$ is called
$M$-structure if $K\in U(\cal L)$ is a central element of the
algebra $U(\cal L)$.

{\bf Lemma 2.3.} For any $M$-structure $\cal L$ we have
$$p^i=-\phi^i_{k,l}\lambda^{l,k},\qquad
q_i=-\psi_i^{k,l}\mu_{l,k},\qquad
u_i^j=-t_i^j-\phi^j_{k,l}\psi_i^{l,k}.$$

{\bf Proof.} This is a direct consequence of the identities
$A_iK=KA_i$ and $B^jK=KB^j$.

{\bf Lemma 2.4.} For $M$-structure $\cal L$ the algebra $U(\cal L)$
is defined by the generators $A_1,...,A_p,$ $B^1,...,B^p$ and
relations obtained from (\ref{(walg)}), (\ref{(walg1)}) by
elimination of $C$. Tensors $\phi_{i,j}^k,$ $\psi_k^{i,j},$ $
\mu_{i,j},$ $\lambda^{i,j}$ should satisfy the properties
(\ref{(!!!!)}), (\ref{(!!!!!)}), (\ref{(!!!!!!)}). Any algebra
defined by such generators and relations is isomorphic to $U(\cal
L)$ for a suitable $M$-structure $\cal L$.

{\bf Theorem 2.2.} Let $\cal L$ be an $M$-structure. Then for any
representation ${U(\cal L)} \to Mat_n$ given by $A_1\to
a_1,...,A_p\to a_p, B^1\to b^1,...,B^p\to b^p,C\to c$ the formula
(\ref{sekprod}) defines an associative product on $Mat_n$ compatible
with the usual product.

{\bf Proof.} Comparing (\ref{(!)})-(\ref{(!!!)}) with
(\ref{(walg)})-(\ref{(walg2)}), where $p^i$, $q_i$ and $u_i^j$ are
given by Lemma 2.3, we see that this is just reformulation of the
Theorem 2.1.

{\bf Definition.} A representation of  $U(\cal L)$ is called
non-degenerate if the matrices $a_1,...,a_p,1$ are linear
independent as well as $b^1,...,b^p,1.$

{\bf Remark.} It is clear that $M$-structure $\cal L^{\prime}$ is
equivalent to $\cal L$ if the defining relations for $U(\cal
L^{\prime})$ can be obtained from the defining relations for
$U(\cal L)$ by  a transformation of the form $$A_i\to g_i^k
A_k+u_i,\quad B^i\to  h_k^i B^k+v^i,\quad C\to C-u_i h_k^i B^k-v^i
g_i^k A_k-u_i v^i$$ where  $u_1,...,u_p, v^1,...,v^p$ are some
constants and $g_i^k h_k^j=\delta_i^j$.

{\bf Theorem 2.3.} There is an one-to-one correspondence between $n$
dimensional non-degenerate representations of algebras $U(\cal L)$
corresponding to $M$-structures up to equivalence of $M$-structures
and associative products on $Mat_n$ compatible with the usual
product.

{\bf Proof.} This is a direct consequence of Theorems 2.1 and 2.2.

The structure of the algebra $U(\cal L)$ for $M$-structure $\cal
L$ is described by the following

{\bf Theorem 2.4.} The algebra $U(\cal L)$ is spanned by the
elements $K^s,\,A_i K^s,\, B_j K^s,\, A_i B^j K^s,$ where
$i,j=1,...,p,\,$ and $\, s=0,1,2,...$

{\bf Proof.} Since $K$ is a central element, we have
$KA_i=A_iK,\,\, KB^j=B^jK,\,\, KC=CK$. Using this, one can check
that a product of any elements listed in the theorem can be
written as a linear combination of these elements. To prove the
theorem one should also check associativity, which is possible to
do directly.

{\bf Remark.} As we have mentioned in the Introduction, if a linear
space $\bf V$ is equipped with two compatible associative
multiplications, then one can construct those in the space
$Mat_m(\bf V)$.  Since $Mat_m(Mat_n)=Mat_{mn}$, in the matrix case
this construction yields a second multiplication for the algebras
$Mat_{mn}, \,\, m=1,2,...$ if we have a second multiplication in
$Mat_n$. One can see that in the language of representations of
$M$-structures this corresponds to the direct sum of $m$ copies of a
given $n$-dimensional representation.

{\bf Example 2.1.} Suppose $\cal A$ and $\cal B$ are generated by
 elements $A\in \cal A$ and  $B\in
\cal B$ such that $A^{p+1}=B^{p+1}=1$.  Take
$1,A,...,A^p,B,...,B^p,C$ for a basis in $\cal L$ and assume that
$(B^i,A^{-i})=\epsilon^i-1$, $(1,C)=1$ and other scalar products are
equal to zero. Here $\epsilon$ is a primitive root of unity of order
$p+1$. The structures of left $\cal B$-module and right $\cal
A$-module on $\cal L$ are defined by the formulas:
$$B^iA^j=\frac{\epsilon^{-j}-1}{\epsilon^{-i-j}-1}A^{i+j}+\frac{\epsilon^i-1}{\epsilon^{i+j}-1}B^{i+j},$$
for $i+j\ne0$ modulo $p$ and
$$B^iA^{-i}=1+(\epsilon^i-1)C,$$
$$CA^i=\frac{1}{1-\epsilon^i}A^i+\frac{1}{\epsilon^i-1}B^i,$$
$$B^iC=\frac{1}{\epsilon^{-i}-1}A^i+\frac{1}{1-\epsilon^{-i}}B^i$$
for $i\ne0$ modulo $p+1$. One can see that these formulas define an
$M$-structure. The central element has the following form
$K=C+\sum_{0<i<p}\frac{1}{\epsilon^i-1}A^{-i}B^i$.

Let $a$, $t$ be linear operators in some vector space. Assume that
$a^{p+1}=1$, $at=\epsilon ta$ and the operator $t-1$ is invertible.
It is easy to check that the formulas $$A\to a,\qquad B\to
\frac{\epsilon t-1}{t-1}a,\qquad C\to\frac{t}{t-1}$$ define a
representation of the algebra $U(\cal L)$. Note that we do not
assume that $t^{p+1}=1$. We have only $at^{p+1}=t^{p+1}a$ which
easily follows from the commutation relation between $a$ and $t$.

\subsection{Case of commutative semi-simple algebra ${\cal A}$}

Consider the case
\begin{equation}\label{AA}
A_i A_j=\delta_{i,j}A_i.
\end{equation}
 In other words
\begin{equation}\label{AA1} \phi^k_{i,j}=\delta_{i,j}\delta_{i,k},\qquad
\mu_{i,j}=0.
\end{equation}

{\bf Theorem 2.5.} In this case any corresponding algebra ${\cal B}$
can be reduced by an appropriate shift $B^{i}\rightarrow
B^{i}+c_{i}$ to one defined by the formulas:
\begin{equation}\label{BB1}
B^i B^j=(u_i-q_{i,j}) B^i+q_{i,j} B^j+v_i, \qquad i\ne j
\end{equation}
and
\begin{equation}\label{BB2}
(B^i)^2=u_i B^i+v_i,
\end{equation}
where constants $u_i,v_i,q_{i,j}$ satisfy the following relations:
\begin{equation}\label{BB3}q_{i,j}^2=u_i q_{i,j}+v_i,
\end{equation}
\begin{equation}\label{BB4}
(u_i-q_{i,j})^2=u_j(u_i-q_{i,j})+v_j,
\end{equation}
where $i\ne j$, and
\begin{equation}\label{BB5}
(q_{i,k}-q_{j,k})(q_{i,k}-q_{i,j})=0
\end{equation}
for pairwise distinct $i,j,k$. The corresponding algebra $U(\cal L)$
is determined by the formulas (\ref{AA}), (\ref{BB1}), (\ref{BB2})
and:
$$B^iA_j=(u_j-q_{j,i})A_j\qquad \hbox{for} \quad i\ne j, \qquad
\quad B^iA_i=u_iA_i+\sum_{k\ne i}q_{k,i}A_k+B^i+C,$$
$$
B^iC=\sum_{1\le k\le p}v_kA_k-u_iB^i-v_i, \qquad CA_j=-u_j A_j.
$$

{\bf Proof.} In our case the first equation of (\ref{(!!!!!!)})
reads $\delta_{j,k}t^i_k=\delta_{i,j}t^i_k,$ which gives
$t^i_j=\delta^i_jr_j$ for some tensor $r_j$. The third equation of
(\ref{(!!!!!)}) reads
$\delta_{j,k}\psi^{l,i}_j=\delta_{l,k}\psi^{l,i}_j+\delta_{i,j}\psi^{l,i}_k+(r_j-
r_k)\delta^l_k\delta^i_j-\delta^i_j\psi^{l,l}_k,$ which has the
following general solution
$\psi^{l,i}_j=\delta^l_j(h_j-r_i-q_{l,i})+\delta^i_jq_{l,i}$.  From
the second equation of (\ref{(!!!!!!)}) we find
$\lambda^{k,j}=\lambda^{k,k}+q_{k,j}(r_j-r_k)$. Substituting these
into the formulas for the product in the algebra $\cal B,$ we get
(\ref{BB1}) and (\ref{BB2}) after suitable shift of $B^1,...,B^p$
for some $u_i, v_i$. Associativity of the algebra $\cal B$ gives
(\ref{BB3}), (\ref{BB4}) and (\ref{BB5}). Indeed, consider an
algebra defined by identities
$$
B_{i} B_{j}=p_{ij} B_{i}+q_{ij} B_{j}+r_{ij}, \quad i\ne j, \qquad
\qquad B_{i}^{2}=u_{i} B_{i}+v_{i}
$$
This algebra is associative iff
$$
r_{ij}=-p_{ij} q_{ij}, \qquad q_{ij}^{2}=u_{i} q_{ij}+v_{i},
$$
$$
(p_{ij}-q_{jk})(p_{ik}-p_{jk})=0, \qquad
(p_{ij}-q_{jk})(q_{ik}-q_{ij})=0,
$$
which equivalent to (\ref{BB3}), (\ref{BB4}) and (\ref{BB5}) in our
case. The explicit form of identities for the algebra $U(\cal L)$
follows from (\ref{(walg1)}) and (\ref{(walg2)}).

{\bf Remark.} It follows from (\ref{BB1}), (\ref{BB2}) that the
vector space spanned by $1$ and $B^i$, where $i$ belongs to
arbitrary subset of the set $\{1,2,\dots,p \},$ is a subalgebra in
$\cal B.$

Two algebras (\ref{BB1}) - (\ref{BB5}) are said to be equivalent if
they are related by a transformation of the form
\begin{equation}\label{trra}
B^{i}\rightarrow c_{1} B^{i}+c_{2}, \qquad i=1,\dots,p
\end{equation}
and a permutation of the generators $B^{1}, \dots, B^{p}.$

{\bf Example 2.2.} Suppose $\cal B$ is commutative. It follows
from
\begin{equation}\label{cBB1}
B^i B^j-B^j B^i=(u_i-q_{i,j}-q_{j,i}) B^i- (u_j-q_{i,j}-q_{j,i})
B^j+(v_i-v_{j})
\end{equation}
that in this case we have $u_{i}=u_{j},$  $v_{i}=v_{j}$ and
$q_{i,j}+q_{j,i}=u_{i}$ for any $i,j$. Such an algebra is equivalent
to the one defined by
$$
u_{1}=\cdots=u_{p}=0, \qquad v_{1}=\cdots=v_{p}=1, \qquad \quad
q_{ij}=1, \quad
 q_{ji}=-1, \quad i>j.
$$
It is easy to verify that this algebra is semi-simple.

{\bf Example 2.3.} One solution of the system (\ref{BB3}) -
(\ref{BB5}) is obvious:
$$
q_{ij}=u_{i}+\tau, \qquad v_{i}=\tau^2+u_{i} \tau, \qquad
i,j=1,\dots, p,
$$
where $\tau$ is arbitrary parameter. Using transformation
(\ref{trra}), we can reduce $\tau$ by zero. Algebra ${\cal B}$
described in this example is called {\it regular}. The corresponding
associative product compatible with the matrix product in
$Mat_{p+1}$ have been independently found by I.Z. Golubchik.

Now our aim is to describe all irregular algebras ${\cal B}$. It
follows from (\ref{BB3}), (\ref{BB4}) that
\begin{equation}\label{qijk}
(q_{ki}-q_{kj}) (q_{ki}+q_{kj}-u_{k})=0
\end{equation}
for any distinct $i,j,k$ and
\begin{equation}\label{uvij}
(q_{ij}-u_{i}) (u_{i}-u_{j})=v_{i}-v_{j}
\end{equation}
for any $i\ne j.$ Formula (\ref{uvij}) implies  that
\begin{equation}\label{qij}
(q_{ij}-q_{ji}-u_{i}+u_{j}) (u_{i}-u_{j})=0.
\end{equation}

We associate with any algebra ${\cal B}$ the following equivalence
relation on the set $\{1,2,\dots,p \}$: $i\thicksim j$ if $u_i=u_j.$
Denote by $m$ the number of equivalence classes. It follows from
(\ref{uvij}) that if $i$ and $j$ belong to the same equivalence
class then $v_i=v_j$. Furthermore, if  $i$ and $j$ belong to
different equivalence classes, we have
\begin{equation}\label{qqij}
q_{ij}=u_i+\frac{v_i-v_j}{u_i-u_j}
\end{equation}
and therefore $q_{ij}$ is well defined function on the set of pairs
of equivalence classes.

Besides $\thicksim,$ we consider one more relation $\thickapprox$
defined as follows: $i \thickapprox j$ if $i=j$ or if $u_i=u_j$ and
$q_{ij}=q_{ji}$. It is easy to derive from (\ref{BB5}) that
$\thickapprox$ is an equivalence relation and the value of $q_{ij}$
does not depend on the choosing of $i,j$ from the equivalence class.

{\bf Case m=1.} Consider the case $m=1$ or, the same $u_{i}=u_{1}$
for any $i$. Denote one of two possible values of $q_{ij}$ by
$-\tau.$ It follows from (\ref{BB3}) that other possible value is
equal to $u_1+\tau$ and $v_1=\tau^2+u_1 \tau$. If ${\cal B}$ is
irregular then $-\tau$ and $u_1+\tau$ are distinct.

Given $u_1, \tau$ any algebra ${\cal B}$ is defined by the following
data: arbitrary clustering of the set $\{1,2,\dots,p \}$ into
equivalence classes $K_1,\dots,K_s$ with respect to $\thickapprox$
and arbitrary function $Q_{ij}$ on the pairs of equivalence classes
with values in $ \{-\tau,u_1+\tau\}$ such that $Q_{\alpha \beta}\ne
Q_{\beta \alpha}$ if $\alpha \ne \beta$. The function $q_{ij}$ is
defined as follows: $q_{i,j}=Q_{ij}$ if $i,j$ belong to different
classes and $q_{i,j}=Q_{\alpha \alpha}$ if $i,j\in K_\alpha$. It can
be verified that the parameters defined as above satisfy (\ref{BB3})
- (\ref{BB5}).

{\bf Case m=2.}  In this case we have two distinct parameters $u_1$
and $u_2$. Let ${\cal K}_1$ and ${\cal K}_2$ be corresponding
equivalence classes with respect to $\thicksim$. Denote
$\frac{v_2-v_1}{u_2-u_1}$ by $\tau$. Then $v_i=\tau^2+u_i \tau.$
Using (\ref{qqij}), we get  $q_{ik}=u_\alpha+\tau$ if $i\in {\cal
K}_{\alpha}$ and $k\nsubseteq {\cal K}_{\alpha}.$ It follows from
(\ref{BB5}) that for any $j,$ $q_{ij}$ may take on values
$u_\alpha+\tau$ or $-\tau$. Suppose $i,j\in {\cal K}_\alpha$ and
$k\nsubseteq {\cal K}_\alpha$; then (\ref{BB5}) yields
$$
(q_{ij}-\tau-u_1)(q_{ij}-\tau-u_2)=0.
$$
Therefore if $\tau+u_1\ne -\tau$ and $\tau+u_2\ne -\tau,$ then
${\cal B}$ is regular. In the case $\tau=-\frac{u_2}{2},$ we have
$q_{ij}=\frac{u_2}{2}$ for any $i\in {\cal K}_2.$ Formula
(\ref{qqij}) implies $q_{ji}=u_1-\frac{u_2}{2}.$ To complete the
description we should define $q_{ij}$ for $i,j\in {\cal K}_1.$ It is
clear that the vector space spanned by $1$ and $B^i$, where $i\in
{\cal K}_1,$ is a subalgebra that belongs to the case $m=1$
described above. It turns out that this subalgebra can be chosen
arbitrarily.

{\bf Case m $\ge$ 3.} In this case all possible algebras $\cal B$
can be described as follows.

{\bf Proposition 2.1.} Suppose $u_{1},\dots,u_{m}$ are  pairwise
distinct and $m\ge 3$. Then if $p > 3$ then $\cal B$ is regular.
For $p=3$ there exists one more algebra described in the following

{\bf Example 2.4.} The algebra $\cal B$ defined by relations
\begin{equation}
\begin{array}{c}
q_{21}=q_{31}, \qquad q_{12}=q_{32}, \qquad q_{13}=q_{23},\\[2mm]
u_{1}=q_{12}+q_{13}, \qquad u_{2}=q_{21}+q_{23},  \qquad
u_{3}=q_{31}+q_{32},\\[2mm]
 v_{1}=-q_{12} q_{13}, \qquad v_{2}=-q_{21} q_{23},  \qquad v_{3}=-q_{31}
q_{32}
\end{array} \label{poss1}
\end{equation}
is isomorphic to $Mat_{2}$.

{\bf Proof.} Suppose that $u_{i},u_{j}, u_{k}$ are pairwise
distinct. Then we deduce from (\ref{qij}) that
$$
q_{ij}-q_{ji}-u_{i}+u_{j}=q_{jk}-q_{kj}-u_{j}+u_{k}=q_{ki}-q_{ik}-u_{k}+u_{i}=0
$$
and therefore $ q_{ij}+q_{jk}+q_{ki}=q_{ji}+q_{ik}+q_{kj}. $ It
follows from this formula and (\ref{BB5}) that there exist only two
possibilities:
\begin{equation}\label{pos1}
 \qquad \qquad q_{ji}=q_{ki}, \qquad q_{ij}=q_{kj}, \qquad q_{ik}=q_{jk}
\end{equation}
or
\begin{equation}\label{pos2}
 \qquad \qquad q_{ij}=q_{ik}, \qquad q_{jk}=q_{ji}, \qquad q_{ki}=q_{kj}.
\end{equation}
It is not difficult to show that  $\cal B$ is regular if it
contains a triple of the type (\ref{pos2}). It can be verified
also that if $\cal B$ contains a triple  of the type (\ref{pos1}),
then $\cal B$ coincides with the algebra described in Example 2.4.

{\bf Remark.} It follows from Theorem 4.2 of Section 4 that if
$\cal A$ is commutative and semi-simple and $\cal B$ is
semi-simple, then either $\cal B$ is commutative (Example 2.2) or
${\cal B}=Mat_2$ (Example 2.4).

\section{Semi-simple case}

Consider an associative algebra $M=\oplus_{1\le\alpha\le m}
M_{\alpha}$, where $M_{\alpha}$ is isomorphic to $Mat_{n_{\alpha}}$.
We are going to study associative products in this algebra
compatible with the usual one. All constructions and results are
similar to the matrix case.

We use Greek letters for indexes related to the direct summands of
$M$. Throughout this section, we keep the following summation
agreement. We sum by repeated Latin indexes and do not sum by
repeated Greek indexes if the opposite is not stated explicitly.
Symbols $\delta_i^j,$ $\delta_{i,j},$ and $\delta^{i,j}$ stand for
the Kronecker delta.

\subsection{Construction of the second product}

Let $R$ be a linear operator in the space $M$. The operator $R$
takes $x_{\alpha}\in M_{\alpha}$ to
$\sum_{\beta}R_{\beta,\alpha}(x_{\alpha}),$ where
$R_{\beta,\alpha}(x_{\alpha})\in M_{\beta}$. It is clear that
$R_{\beta,\alpha}$ is a linear map from $M_{\alpha}$ to
$M_{\beta}$. Note that any linear map from the space
$Mat_{n_{\alpha}}$ to the space $Mat_{n_{\beta}}$ can be written
in the form $x_{\alpha}\to
a_{i,\beta,\alpha}x_{\alpha}b^i_{\alpha,\beta},$ where
$a_{1,\beta,\alpha},...,a_{l,\beta,\alpha}$ are some
$n_{\beta}\times n_{\alpha}$ matrices and
$b^1_{\alpha,\beta},...,b^l_{\alpha,\beta}$ are some
$n_{\alpha}\times n_{\beta}$ matrices.

Assume that
$R_{\beta,\alpha}(x_{\alpha})=a_{i,\beta,\alpha}x_{\alpha}b^i_{\alpha,\beta}$
for $\alpha\ne\beta$ and
$R_{\alpha,\alpha}(x_{\alpha})=a_{i,\alpha,\alpha}x_{\alpha}b^i_{\alpha,\alpha}+c_{\alpha}x_{\alpha}$
for some matrices $a_{i,\beta,\alpha},b^i_{\alpha,\beta},1\le i\le
p_{\alpha,\beta}$ and $c_{\alpha}$, $1\le\alpha,\beta\le m$ with
$p_{\alpha,\beta}$ being smallest possible in the equivalence class
 of $R$. This means that the following sets of
matrices are linear independent:
$\{a_{1,\beta,\alpha},...,a_{p_{\alpha,\beta},\beta,\alpha}\}$,
$\{b^1_{\alpha,\beta},...,b^{p_{\alpha,\beta}}_{\alpha,\beta}\}$
for $\alpha\ne\beta$ and
$\{1,a_{1,\alpha,\alpha},...,a_{p_{\alpha,\alpha},\alpha,\alpha}\}$,
$\{1,b^1_{\alpha,\alpha},...,b^{p_{\alpha,\alpha}}_{\alpha,\alpha}\}$.
It follows from (\ref{mult2}) that  the second product of
$x_{\alpha}\in M_{\alpha}$ and $y_{\beta}\in M_{\beta}$ has the
form
\begin{equation}\begin{array}{c}
\label{sekprod1} x_{\alpha}\circ y_{\beta}=a_{i,\beta,\alpha}
x_{\alpha} b^i_{\alpha,\beta} y_{\beta}+x_{\alpha}
a_{i,\alpha,\beta} y_{\beta}
b^i_{\beta,\alpha}, \qquad \alpha\ne\beta, \\[4mm]
x_{\alpha}\circ y_{\alpha}=a_{i,\alpha,\alpha} x_{\alpha}
b^i_{\alpha,\alpha} y_{\alpha}+x_{\alpha} a_{i,\alpha,\alpha}
y_{\alpha} b^i_{\alpha,\alpha}-a_{i,\alpha,\alpha} x_{\alpha}
y_{\alpha} b^i_{\alpha,\alpha}+x_{\alpha} c_{\alpha} y_{\alpha}.
\end{array}
\end{equation}

We have the following transformations preserving  the equivalence
class of $R$. The first family of such transformations is defined
by
$$a_{i,\alpha,\alpha}\to a_{i,\alpha,\alpha}+u_{i,\alpha}, \qquad
b^i_{\alpha,\alpha}\to b^i_{\alpha,\alpha}+ v^i_{\alpha}, \qquad
c_{\alpha}\to
c_{\alpha}-u_{i,\alpha}b^i_{\alpha,\alpha}-v^i_{\alpha}a_{i,\alpha,\alpha}-u_{i,\alpha}v^i_{\alpha}
$$
for any constants
$u_{1,\alpha},...,u_{p_{\alpha,\alpha},\alpha},v^1_{\alpha},...,v^{p_{\alpha,\alpha}}_{\alpha}$
and the second one is given by
$$
a_{i,\beta,\alpha}\to g_{i,\alpha,\beta}^ka_{k,\beta,\alpha}, \qquad
b^i_{\alpha,\beta}\to h_{k,\alpha,\beta}^ib^k_{\alpha,\beta}, \qquad
c_{\alpha}\to c_{\alpha},
$$
 where
$g_{i,\alpha,\beta}^kh_{k,\alpha,\beta}^j=\delta_i^j$. This means
that we can regard $a_{i,\beta,\alpha}$ and $b^i_{\alpha,\beta}$ as
bases in dual vector spaces.

{\bf Theorem 3.1. } If $\circ$ is an associative product on the
space $M,$ then
\begin{equation}\label{(prod1)}
a_{i,\alpha,\beta}a_{j,\beta,\gamma}=\phi^k_{i,j,\alpha,\beta,\gamma}a_{k,\alpha,\gamma}+
\delta_{\alpha,\gamma}\mu_{i,j,\alpha,\beta},\qquad
b^i_{\alpha,\beta}b^j_{\beta,\gamma}=\psi^{i,j}_{k,\alpha,\beta,\gamma}b^k_{\alpha,\gamma}+
\delta_{\alpha,\gamma}\lambda^{i,j}_{\alpha,\beta},
\end{equation}
\begin{equation}\label{(prod2)}
b^i_{\alpha,\beta}a_{j,\beta,\gamma}=\phi^i_{j,k,\beta,\gamma,\alpha}b^k_{\alpha,\gamma}+
\psi^{k,i}_{j,\gamma,\alpha,\beta}a_{k,\alpha,\gamma}+\delta_{\alpha,\gamma}t^i_{j,\alpha,\beta}+
\delta_{\alpha,\gamma}\delta^i_jc_{\alpha},
\end{equation}
\begin{equation}\label{(prod3)}
c_{\alpha}a_{i,\alpha,\beta}=\mu_{i,k,\alpha,\beta}b^k_{\alpha,\beta}-t^k_{i,\beta,\alpha}a_{k,\alpha,\beta}
-\sum_{1\le\nu\le
m}\phi^k_{l,s,\alpha,\nu,\beta}\psi^{s,l}_{i,\beta,\nu,\alpha}a_{k,\alpha,\beta}-
\sum_{1\le\nu\le
m}\delta_{\alpha,\beta}\mu_{l,s,\alpha,\nu}\psi^{s,l}_{i,\alpha,\nu,\alpha},
\end{equation}
\begin{equation}\label{(prod4)}
b^i_{\alpha,\beta}c_{\beta}=\lambda^{k,i}_{\beta,\alpha}a_{k,\alpha,\beta}-
t^i_{k,\alpha,\beta}b^k_{\alpha,\beta}-\sum_{1\le\nu\le m}
\phi^i_{l,s,\beta,\nu,\alpha}\psi^{s,l}_{k,\alpha,\nu,\beta}b^k_{\alpha,\beta}-\sum_{1\le\nu\le
m}
\delta_{\alpha,\beta}\phi^i_{s,l,\alpha,\nu,\alpha}\lambda^{l,s}_{\alpha,\nu,\alpha},
\end{equation}
where $\phi^k_{i,j,\alpha,\beta,\gamma}, \mu_{i,j,\alpha,\beta},
\psi^{i,j}_{k,\alpha,\beta,\gamma}, \lambda^{i,j}_{\alpha,\beta},
t^i_{j,\alpha,\beta}$ are tensors satisfying the properties
\begin{equation}\label{(eq1)}
\phi^s_{j,k,\alpha,\beta,\gamma}\phi^i_{s,l,\alpha,\gamma,\delta}+
\delta^i_l\delta_{\alpha,\gamma}\mu_{j,k,\alpha,\beta}=\phi^i_{j,s,\alpha,\beta,\delta}\phi^s_{k,l\beta,\gamma,\delta}
+\delta^i_j\delta_{\beta,\delta}\mu_{k,l,\beta,\gamma},
 \end{equation}

 $$\phi^s_{j,k\alpha,\beta,\gamma}\mu_{i,s,\alpha,\gamma}=
 \phi^s_{i,j,\beta,\gamma,\alpha}\mu_{s,k,\alpha,\beta},$$

 \begin{equation}\label{(eq2)}
 \psi^{i,j}_{s,\alpha,\beta,\gamma}\psi^{s,k}_{l,\alpha,\gamma,\delta}+
 \delta^k_l\delta_{\alpha,\gamma}\lambda^{i,j}_{\alpha,\beta}
 =\psi^{j,k}_{s,\beta,\gamma,\delta}\psi^{i,s}_{l,\alpha,\beta,\delta}+
 \delta^i_l\delta_{\beta,\delta}\lambda^{j,k}_{\beta,\gamma},
 \end{equation}

 $$
 \psi_{s,\alpha,\beta,\gamma}^{i,j}\lambda^{s,k}_{\alpha,\gamma}=
 \psi^{j,k}_{s,\beta,\gamma,\alpha}\lambda^{i,s}_{\alpha,\beta},$$

 $$\begin{array}{c}
 \phi^s_{j,k,\beta,\gamma,\delta}\psi_{s,\delta,\alpha,\beta}^{l,i}=
 \phi^l_{s,k,\alpha,\gamma,\delta}\psi_{j,\gamma,\alpha,\beta}^{s,i}+
 \phi^i_{j,s,\beta,\gamma,\alpha}\psi_{k,\delta,\alpha,\gamma}^{l,s}+\\ [4mm]
 \delta^l_k\delta_{\alpha,\gamma}t^i_{j,\alpha,\beta}-
 \delta^i_j\delta_{\alpha,\gamma}t_{k,\delta,\alpha}^l-
 \delta^i_j\delta_{\alpha,\gamma}\sum_{1\le\nu\le
m}\phi^l_{s,r,\alpha,\nu,\delta}\psi^{r,s}_{k,\delta,\nu,\alpha},
 \end{array}
$$

 \begin{equation}\label{(eq3)}
 \phi^s_{j,k,\beta,\gamma,\alpha}t^i_{s,\alpha,\beta}=
 \psi_{j,\gamma,\alpha,\beta}^{s,i}\mu_{s,k,\alpha,\gamma}+
 \phi^i_{j,s,\beta,\gamma,\alpha}t_{k,\alpha,\gamma}^s-
 \delta^i_j\delta_{\alpha,\gamma}\sum_{1\le\nu\le
m}\psi_{k,\alpha,\nu,\alpha}^{s,r}\mu_{r,s,\alpha,\nu},
 \end{equation}

 $$\psi^{k,i}_{s,\alpha,\beta,\gamma}t_{j,\alpha,\gamma}^s=
 \phi^i_{j,s,\gamma,\alpha,\beta}\lambda^{k,s}_{\alpha,\beta}+
 \psi_{j,\alpha,\beta,\gamma}^{s,i}t_{s,\alpha,\beta}^k-
 \delta^i_j\delta_{\alpha,\beta}\sum_{1\le\nu\le
m}\phi^k_{s,r,\alpha,\nu,\alpha}\lambda^{r,s}_{\alpha,\nu}$$

{\bf Proof} is similar to the matrix case. Instead of Lemma 2.1
one can use the following

{\bf Lemma 3.1.} Let $x\to p_1xq_1+...+p_lxq_l$ be a zero map from
$Mat_{\alpha}$ to $Mat_{\beta}$. If $p_1,...,p_l$ are linear
independent matrices, then $q_1=...=q_l=0$. Similarly, if
$q_1,...,q_l$ are linear independent matrices, then
$p_1=...=p_l=0$.

\subsection{$PM$-structures and corresponding associative algebras}

In this subsection we describe the algebraic structure underlying
Theorem 3.1.

{\bf Definition.} By weak $PM$-structure (of size $m$) on a linear
space $\cal L$ we mean the following data.
\begin{itemize}
\item  Two subspaces $\cal A$ and $\cal B$ and a distinguished
element $1\in\cal A\cap\cal B\subset \cal L$.

\item  A non-degenerate symmetric scalar product $(\cdot, \cdot)$ on the
space $\cal L$.

\item  Associative products $\cal A\times\cal A\to\cal A$ and
$\cal B\times\cal B\to\cal B$ with unity $1$.

\item  A left action $\cal B\times\cal L\to\cal L$ of the algebra
$\cal B$ and a right action $\cal L\times\cal A\to\cal L$ of the
algebra $\cal A$ on the space $\cal L,$ which commute with each
other.
\end{itemize}

These data should satisfy the following properties:

{\bf 1.} $\dim{\cal A\cap\cal B=\dim\cal L/(\cal A+\cal B)}= m$.
The intersection of $\cal A$ and $\cal B$ is a $m$-dimensional
algebra isomorphic to $\C\oplus...\oplus\C$.

{\bf 2.} The restriction of the action $\cal B\times\cal L\to\cal L$
to the subspace $\cal B\subset \cal L$ is the product in $\cal B$.
The restriction of the action $\cal L\times\cal A\to\cal L$ to the
subspace $\cal A\subset \cal L$ is the product in $\cal A$.

{\bf 3.} $(a_1,a_2)=(b_1,b_2)=0$ and $(b_1b_2,v)=(b_1,b_2v)$,
$(v,a_1a_2)=(va_1,a_2)$ for any $a_1,a_2\in\cal A$, $b_1,b_2\in\cal
B$ and $v\in\cal L$.

It follows from these properties that $(\cdot,\cdot)$ defines a
non- degenerate pairing between ${\cal A}/\cal A\cap\cal B$ and
${\cal B}/\cal A\cap\cal B$, so $\dim\cal A=\dim\cal B$ and
$\dim{\cal L}=2\dim{\cal A}$.

{\bf Lemma 3.2.} Let $\{e_{\alpha};1\le\alpha\le m\}$ be a basis of
the space $\cal A\cap\cal B$ such that
\begin{equation}\label{(un)}
e_{\alpha}e_{\beta}=\delta_{\alpha,\beta}e_{\alpha}.
\end{equation}
Denote by ${\cal L}_{\alpha,\beta}$ the vector space consisting of
elements $v_{\alpha,\beta} \in\cal L$ with the property
\begin{equation}\label{(un1)}
e_{\alpha} v_{\alpha,\beta}=v_{\alpha,\beta}
e_{\beta}=v_{\alpha,\beta}.
\end{equation}
Let $\cal A_{\alpha,\beta}=\cal A\cap\cal L_{\alpha,\beta}$ and
$\cal B_{\alpha,\beta}=\cal B\cap\cal L_{\alpha,\beta}$. Then the
following properties hold:
\begin{itemize}
\item   ${\cal L}=\oplus_{1\le\alpha,\beta\le m}\cal
L_{\alpha,\beta}$,\quad ${\cal A}=\oplus_{1\le\alpha,\beta\le
m}\cal A_{\alpha,\beta}$ and ${\cal B}=\oplus_{1\le\alpha,\beta\le
m}\cal B_{\alpha,\beta}.$
\item  $\dim{\cal A_{\alpha,\beta}\cap\cal
B_{\alpha,\beta}=\dim\cal L/(\cal A_{\alpha,\beta}+\cal
B_{\alpha,\beta})}= \delta_{\alpha,\beta}$. The intersection of
$\cal A_{\alpha,\alpha}$ and $\cal B_{\alpha,\alpha}$ is an
one-dimensional space spanned by $e_{\alpha}$.
\item ${\cal B}_{\alpha,\beta}{\cal
L}_{\beta^{\prime},\gamma}=0$ for $\beta\ne\beta^{\prime}$ and $\cal
B_{\alpha,\beta}\cal L_{\beta,\gamma}\subset \cal
L_{\alpha,\gamma}$. Similarly ${\cal L}_{\alpha,\beta}{\cal
A}_{\beta^{\prime},\gamma}=0$ for $\beta\ne\beta^{\prime}$ and $\cal
L_{\alpha,\beta}\cal A_{\beta,\gamma}\subset \cal
L_{\alpha,\gamma}$. In particular, ${\cal A}_{\alpha,\beta}{\cal
A}_{\beta^{\prime},\gamma}={\cal B}_{\alpha,\beta}{\cal
B}_{\beta^{\prime},\gamma}=0$ for $\beta\ne\beta^{\prime}$ and $\cal
A_{\alpha,\beta}\cal A_{\beta,\gamma}\subset \cal
A_{\alpha,\gamma}$, $\cal B_{\alpha,\beta}\cal
B_{\beta,\gamma}\subset \cal B_{\alpha,\gamma}$.
\item ${\cal L}_{\alpha,\beta}\bot {\cal
L}_{\beta^{\prime},\alpha^{\prime}}$ if $\alpha\ne\alpha^{\prime}$
or $\beta\ne\beta^{\prime}$.
\end{itemize}

It follows from these properties that $(\cdot,\cdot)$ gives non-
degenerate pairing between $\cal A_{\alpha,\beta}$ and $\cal
B_{\beta,\alpha}$ for $\alpha\ne\beta$ and between ${\cal
A}_{\alpha,\alpha}/\C e_{\alpha}$ and ${\cal B}_{\alpha,\alpha}/\C
e_{\alpha}$. Therefore $\dim\cal A_{\alpha,\beta}=\dim\cal
B_{\beta,\alpha}.$

{\bf Proof.} It is clear that $1=e_1+...+e_m$. For $v\in\cal L$ we
have $v=1v1=\sum_{1\le\alpha,\beta\le m}e_{\alpha}ve_{\beta}$ and
$e_{\alpha}ve_{\beta}\in\cal L_{\alpha,\beta},$ which proves all
statements of Lemma 3.2.

{\bf Definition.} By weak $PM$-algebra associated with a weak
$PM$-structure $\cal L$ we mean an associative algebra $U(\cal L)$
possessing the following properties:

{\bf 1.} $\cal L\subset U(\cal L)$ and the actions $\cal
B\times\cal L\to\cal L$, $\cal L\times\cal A\to\cal L$ are the
restrictions of the product in $U(\cal L)$.

{\bf 2.} For any algebra $X$ with the property {\bf 1} there
exists a unique homomorphism of algebras $X\to U(\cal L)$
identical on $\cal L$.

It is easy to see that if $U(\cal L)$ exists, then it is unique for
given $\cal L$. Let us describe the structure of $U(\cal L)$
explicitly. Let $\{e_{\alpha},A_{i,\alpha,\alpha}; 1\le i\le
p_{\alpha,\alpha}\}$ be a basis of $\cal A_{\alpha,\alpha}$ and
$\{e_{\alpha},B^i_{\alpha,\alpha}; 1\le i\le p_{\alpha,\alpha}\}$ be
the dual basis of $\cal B_{\alpha,\alpha}$. Let
$\{A_{i,\alpha,\beta}; 1\le i\le p_{\beta,\alpha}\}$ be a basis of
$A_{\alpha,\beta}$ for $\alpha\ne\beta$ and $\{B^i_{\beta,\alpha};
1\le i\le p_{\beta,\alpha}\}$ be the dual basis of
$B_{\beta,\alpha}$. This means that
$(A_{i,\alpha,\beta},B^j_{\beta^{\prime},\alpha^{\prime}})=
\delta^j_i\delta_{\alpha,\alpha^{\prime}}\delta_{\beta,\beta^{\prime}}$.
Take $C_{\alpha}\in\cal L_{\alpha,\alpha}$ that does not belong to
the sum of $\cal A_{\alpha,\alpha}$ and $\cal B_{\alpha,\alpha}$.
Since $(\cdot , \cdot)$ is non-degenerate, we have
$(e_{\alpha},C_{\alpha})\ne 0$. Multiplying $C_{\alpha}$ by
constant, we can assume that $(e_{\alpha},C_{\alpha})=1$. Adding a
linear combination of
$e_{\alpha},A_{1,\alpha,\alpha},...,A_{p_{\alpha,\alpha},\alpha,\alpha},
B^1_{\alpha,\alpha},...,B_{\alpha,\alpha}^{p_{\alpha,\alpha}}$ to
$C_{\alpha}$, we can assume that
$(C_{\alpha},C_{\alpha})=(C_{\alpha},A_{i,\alpha,\alpha})=(C_{\alpha},B^j_{\alpha,\alpha})=0.$
Such element $C_{\alpha}$ is uniquely determined by choosing of
basis in $\cal A_{\alpha,\alpha}.$

{\bf Lemma 3.3.} The algebra $U(\cal L)$ is defined by (\ref{(un)}),
(\ref{(un1)}) and the following relations
\begin{equation}\label{(rel1)}
A_{i,\alpha,\beta}A_{j,\beta,\gamma}=\phi^k_{i,j,\alpha,\beta,\gamma}A_{k,\alpha,\gamma}+
\delta_{\alpha,\gamma}\mu_{i,j,\alpha,\beta},\qquad
B^i_{\alpha,\beta}B^j_{\beta,\gamma}=\psi^{i,j}_{k,\alpha,\beta,\gamma}B^k_{\alpha,\gamma}+
\delta_{\alpha,\gamma}\lambda^{i,j}_{\alpha,\beta}
\end{equation}

\begin{equation}\label{(rel2)}
B^i_{\alpha,\beta}A_{j,\beta,\gamma}=\phi^i_{j,k,\beta,\gamma,\alpha}B^k_{\alpha,\gamma}+
\psi^{k,i}_{j,\gamma,\alpha,\beta}A_{k,\alpha,\gamma}+\delta_{\alpha,\gamma}t^i_{j,\alpha,\beta}+
\delta_{\alpha,\gamma}\delta^i_jC_{\alpha}
\end{equation}

\begin{equation}\label{(rel3)}
C_{\alpha}A_{i,\alpha,\beta}=\mu_{i,k,\alpha,\beta}B^k_{\alpha,\beta}+u^k_{i,\beta,\alpha}A_{k,\alpha,\beta}
+\delta_{\alpha,\beta}p^i_{\alpha}
\end{equation}

\begin{equation}\label{(rel4)}
B^i_{\alpha,\beta}C_{\beta}=\lambda^{k,i}_{\beta,\alpha}A_{k,\alpha,\beta}+u^i_{k,\alpha,\beta}B^k_{\alpha,\beta}
+\delta_{\alpha,\beta}q_{i,\alpha}
\end{equation}
for certain tensors $\phi^k_{i,j,\alpha,\beta,\gamma},
\mu_{i,j,\alpha,\beta}, \psi^{i,j}_{k,\alpha,\beta,\gamma},
\lambda^{i,j}_{\alpha,\beta},
t^i_{j,\alpha,\beta},u^k_{i,\beta,\alpha},p^i_{\alpha},q_{i,\alpha}$.

{\bf Proof.} Relations (\ref{(rel1)}) just mean that $\cal A$ and
$\cal B$ are associative algebras. Since $\cal L$ is a left $\cal
B$-module and a right $\cal A$-module,
$B^i_{\alpha,\beta}A_{j,\beta,\gamma},C_{\alpha}A_{j,\alpha,\beta},B^i_{\alpha,\beta}C_{\beta}$
should be linear combinations of the basis elements in $\cal L$.
Applying properties {\bf 1, 2, 3} of weak $PM$-structure and Lemma
3.2, we obtain required form of these products. The universality
condition of $U(\cal L)$ shows that this algebra is defined by
(\ref{(un)}), (\ref{(un1)}), (\ref{(rel1)}) - (\ref{(rel4)}).

It is clear that $U({\cal L})=\oplus_{1\le\alpha,\beta\le m}U(\cal
L)_{\alpha,\beta}$, where $U({\cal L})_{\alpha,\beta}=\{v\in
U({\cal L}); e_{\alpha}v=ve_{\beta}=v\}$. We have $U({\cal
L})_{\alpha,\beta}U({\cal L})_{\beta^{\prime},\gamma}=0$ for
$\beta\ne\beta^{\prime}$ and $U({\cal L})_{\alpha,\beta}U({\cal
L})_{\beta,\gamma}\subset U({\cal L})_{\alpha,\gamma}$.

Let us define an element $K_{\alpha}\in U(\cal L)$ by the formula
$K_{\alpha}=C_{\alpha}+\sum_{1\le\nu\le m}A_{i,\alpha,\nu}
B_{\nu,\alpha}^i$. It is clear that $K_{\alpha}$ thus defined does
not depend on the choice of the basis in $\cal A$ and $\cal B$
provided $(A_{i,\alpha,\beta},B^j_{\beta^{\prime},\alpha^{\prime}})=
\delta^j_i\delta_{\alpha,\alpha^{\prime}}\delta_{\beta,\beta^{\prime}}$,
$(e_{\alpha},C_{\alpha})=1$ and
$(C_{\alpha},C_{\alpha})=(C_{\alpha},A_{i,\alpha,\alpha})=(C_{\alpha},B^j_{\alpha,\alpha})=0$.
Indeed, the coefficients of $K_{\alpha}$ are just entries of the
tensor inverse to the form $(\cdot, \cdot)$.

{\bf Definition.} A weak $PM$-structure $\cal L$ is called
$PM$-structure if $K=\sum_{1\le\alpha\le m}K_{\alpha}\in U(\cal L)$
is a central element of the algebra $U(\cal L)$.

It is clear that $K$ is central if and only if
$K_{\alpha}v=vK_{\beta}$ for all $v\in U({\cal L})_{\alpha,\beta}$.

{\bf Lemma 3.4.} For any $PM$-structure $\cal L,$ we have
$$p^i_{\alpha}=-\sum_{1\le\nu\le m}
\phi^i_{s,l,\alpha,\nu,\alpha}\lambda^{l,s}_{\alpha,\nu,\alpha},\qquad
q_{i,\alpha}=- \sum_{1\le\nu\le
m}\mu_{l,s,\alpha,\nu}\psi^{s,l}_{i,\alpha,\nu,\alpha},$$
$$
u_{i,\alpha,\beta}^j=- t^j_{i,\alpha,\beta}-\sum_{1\le\nu\le m}
\phi^j_{l,s,\beta,\nu,\alpha}\psi^{s,l}_{i,\alpha,\nu,\beta}$$

{\bf Proof.} This is a direct consequence of
$A_{i,\alpha,\beta}K_{\beta}=K_{\alpha}A_{i,\alpha,\beta}$ and
$B^j_{\alpha,\beta}K_{\beta}=K_{\alpha}B^j_{\alpha,\beta}$.

{\bf Lemma 3.5.} For any $PM$-structure $\cal L,$ the algebra
$U(\cal L)$ is defined by the generators \linebreak
$\{e_{\alpha},A_{i,\alpha,\beta},B^i_{\beta,\alpha}; 1\le i\le
p_{\beta,\alpha},1\le\alpha,\beta\le m\}$ and relations obtained
from (\ref{(un)}), (\ref{(rel1)}), (\ref{(rel2)}) by elimination of
$C_{\alpha}$. Tensors
$\phi_{i,j}^k,\psi_k^{i,j},\mu_{i,j},\lambda^{i,j}$ should satisfy
the properties (\ref{(eq1)})-(\ref{(eq3)}). Any algebra defined by
such generators and relations is isomorphic to $U(\cal L)$ for a
suitable $PM$-structure $\cal L$.

Let $\tau : U({\cal L})\to End(V)$ be a representation of the
algebra $U(\cal L)$. Let $\pi_{\alpha}=\tau(e_{\alpha})$ and
$V_{\alpha}=\pi_{\alpha}(V)$. It follows from (\ref{(un)}) that
$V=\oplus_{1\le\alpha\le m}V_{\alpha}$. Let
$n_{\alpha}=dim(V_{\alpha})$. We can regard $x\in U(\cal
L)_{\alpha,\beta}$ as a linear operator from $V_{\beta}$ to
$V_{\alpha}$ or, choosing basis in $V_1,...,V_m$, as an
$n_{\alpha} \times n_{\beta}$ matrix.

{\bf Definition.} By a representation of a $PM$-algebra $U(\cal L)$
of dimension $(n_1,...,n_m)$ we mean a correspondence
$A_{i,\beta,\alpha}\to a_{i,\beta,\alpha},\,B^i_{\alpha,\beta}\to
b^i_{\alpha,\beta},\, C_{\alpha}\to c_{\alpha};\,\, 1\le i\le
p_{\alpha,\beta},1\le\alpha,\beta\le m,$ where
$a_{i,\alpha,\beta},\,b^i_{\alpha,\beta}$ are $n_{\alpha}\times
n_{\beta}$ matrices and $c_{\alpha}$ are $n_{\alpha}\times
n_{\alpha}$ matrices satisfying (\ref{(prod1)}), (\ref{(prod2)}),
(\ref{(prod3)}), (\ref{(prod4)}).

It is clear  that this definition is equivalent to the usual one for
the associative algebra $U(\cal L)$.

{\bf Theorem 3.2.} Let $\cal L$ be a $PM$-structure. Then for any
representation of $U(\cal L)$ given by $A_{i,\beta,\alpha}\to
a_{i,\beta,\alpha},B^i_{\alpha,\beta}\to b^i_{\alpha,\beta},
C_{\alpha}\to c_{\alpha}; \,\,1\le i\le
p_{\alpha,\beta},1\le\alpha,\beta\le m\,$ the formula
(\ref{sekprod1}) defines an associative product on
$M=\oplus_{1\le\alpha\le m} Mat_{n_{\alpha}}$ compatible with the
usual one.

{\bf Proof.} Comparing (\ref{(prod1)})-(\ref{(prod4)}) with
(\ref{(rel1)})-(\ref{(rel4)}), where $p^i_{\alpha},q_{i,\alpha}$
and $u_{i,\alpha,\beta}^j$ are given by Lemma 3.4, we see that
this is just reformulation of Theorem 3.1.

{\bf Definition.} A representation of a $PM$-algebra $U(\cal L)$ is
called non-degenerate if the following sets of matrices are linear
independent: $\{1,a_{i,\alpha,\alpha}; 1\le i\le
p_{\alpha,\alpha}\}$, $\{1,b_{i,\alpha,\alpha}; 1\le i\le
p_{\alpha,\alpha}\}$, $\{a_{i,\beta,\alpha}; 1\le i\le
p_{\alpha,\beta}\}$ and $\{b ^i_{\alpha,\beta}; 1\le i\le
p_{\alpha,\beta}\}$ for $\alpha\ne\beta$.

{\bf Theorem 3.3. } There is a one-to-one correspondence between
$(n_1,...,n_m)$-dimensional non-degenerate representations of
$PM$-algebras $U(\cal L)$ up to equivalence of the algebras and
associative products on $Mat_{n_1}\oplus...\oplus Mat_{n_m}$
compatible with the usual product.

{\bf Proof.} This is a direct consequence of Theorems 3.1 and 3.2.

The structure of a $PM$-algebra $U(\cal L)$ can be described as
follows.

{\bf Theorem 3.4.} A basis of $U(\cal L)_{\alpha,\beta}$ for
$\alpha\ne\beta$ consists of the elements
$$\{A_{i,\alpha,\beta}K_{\beta}^s,\,\,\,
B^j_{\alpha,\beta}K_{\beta}^s,\,\,\,
A_{i_1,\alpha,\nu}B^{j_1}_{\nu,\beta}K_{\beta}^s\},$$ where $ 1\le
i\le p_{\beta,\alpha},\,\, 1\le j\le p_{\alpha,\beta},\,\,
1\le\alpha,\beta,\nu\le m, 1\le i_1\le p_{\nu,\alpha},\,\, 1\le
j_1\le p_{\nu,\beta},\,\, s=0,1,2,...\,$. A basis of $U(\cal
L)_{\alpha,\alpha}$ consists of the elements $$\{e_{\alpha},\,\,\,
A_{i,\alpha,\alpha}K_{\alpha}^s, \,\,\,
B^j_{\alpha,\alpha}K_{\alpha}^s,\,\,\,
A_{i_1,\alpha,\nu}B^{j_1}_{\nu,\alpha}K_{\alpha}^s\},$$ where $1\le
i,j\le p_{\alpha,\alpha},\,\, 1\le\nu\le m,\,\, 1\le i_1,j_1\le
p_{\nu,\alpha},\,\, s=0,1,2,...\,$.

{\bf Proof.} Since $K$ is a central element, we have
$K_{\alpha}A_{i,\alpha,\beta}=A_{i,\alpha,\beta}K_{\beta},\quad$ $
K_{\alpha}B^j_{\alpha,\beta}=B^j_{\alpha,\beta}K_{\beta},
$\linebreak $ K_{\alpha}C_{\alpha}=C_{\alpha}K_{\alpha}$. Using
this, one can check that a product of any elements listed in the
theorem can be written as a linear combination of these elements. To
prove the theorem, one should also check the associativity, which is
possible to do directly.

{\bf Definition.} Let ${\cal L}_1$ and ${\cal L}_2$ be weak
$PM$-structures. Let ${\cal A}_1,{\cal B}_1\subset{\cal L}_1$ and
${\cal A}_2,{\cal B}_2\subset{\cal L}_2$ be corresponding algebras
and $(\cdot,\cdot)_1$, $(\cdot,\cdot)_2$ be corresponding scalar
products. By direct sum of ${\cal L}_1$ and ${\cal L}_2$ we mean the
weak $PM$-structure ${\cal L}={\cal L}_1\oplus {\cal L}_2$ with
${\cal A}={\cal A}_1\oplus {\cal A}_2$, ${\cal B}={\cal B}_1\oplus
{\cal B}_2$ and $(\cdot,\cdot)=(\cdot,\cdot)_1+(\cdot,\cdot)_2$. We
assume the componentwise action of $\cal A$ and $\cal B$ on $\cal
L$.

{\bf Definition.} A weak $PM$-structure is called indecomposable if
it is not equal to ${\cal L}_1\oplus {\cal L}_2$ for nonzero ${\cal
L}_1$ and ${\cal L}_2$.

It is clear that decomposable $PM$-structures correspond to
decomposable pairs of compatible associative products.

{\bf Definition.} Let ${\cal L}$ be a weak $PM$-structure. By the
opposite weak $PM$-structure ${\cal L}^{op}$ we mean a
$PM$-structure with the same linear space ${\cal L}$, the same
scalar product and algebras $\cal A$, $\cal B$ replaced by the
opposite algebras ${\cal B}^{op}$, ${\cal A}^{op},$ correspondingly.
We remind that a right module over an associative algebra is left
module over opposite algebra and vice-versa.

Let us describe the $PM$-structure related to Example 1.3.

{\bf Example 3.1.} Let $\dim {\cal A}_{\alpha,\beta}=\dim {\cal
B}_{\alpha,\beta}=1$ for all $1\le\alpha,\beta\le m$. Suppose that
for $\alpha\ne\beta$ the space $\cal A_{\alpha,\beta}$ is spanned by
an element $A_{\alpha,\beta}$ and the space $\cal B_{\alpha,\beta}$
is spanned by an element $B_{\alpha,\beta}$. Note that ${\cal
A}_{\alpha,\alpha}={\cal B}_{\alpha,\alpha}=\C e_{\alpha}$. Assume
that $A_{\alpha,\beta}A_{\beta,\gamma}=A_{\alpha,\gamma}$,
$B_{\alpha,\beta}B_{\beta,\gamma}=B_{\alpha,\gamma}$ for
$\alpha\ne\gamma$ and
$A_{\alpha,\beta}A_{\beta,\alpha}=B_{\alpha,\beta}B_{\beta,\alpha}=e_{\alpha}$.
Note that $\dim {\cal L}_{\alpha,\beta}=2$ for all
$1\le\alpha,\beta\le m$ and a basis of $\cal L_{\alpha,\beta}$ is
$\{A_{\alpha,\beta},B_{\alpha,\beta}\}$ for $\alpha\ne\beta$. A
basis of $\cal L_{\alpha,\alpha}$ is $\{e_{\alpha},C_{\alpha}\}$.
Assume that
$(A_{\alpha,\beta},B_{\beta,\alpha})=(u_{\alpha}-u_{\beta})/t_{\beta}$,
$(e_{\alpha},C_{\alpha})=t_{\alpha}^{-1}$ and structures of left
$\cal B$-module and right $\cal A$-module are given by the formulas:
$$B_{\alpha,\beta}A_{\beta,\gamma}=\frac{u_{\beta}-u_{\gamma}}{u_{\alpha}-u_{\gamma}}A_{\alpha,\gamma}
+\frac{u_{\beta}-u_{\alpha}}{u_{\gamma}-u_{\alpha}}B_{\alpha,\gamma}$$
for $\alpha\ne\gamma$ and
$$B_{\alpha,\beta}A_{\beta,\alpha}=e_{\alpha}+(u_{\beta}-u_{\alpha})C_{\alpha},$$
$$C_{\alpha}A_{\alpha,\beta}=\frac{1}{u_{\alpha}-u_{\beta}}A_{\alpha,\beta}+\frac{1}{u_{\beta}
-u_{\alpha}}B_{\alpha,\beta},\qquad
B_{\alpha,\beta}C_{\beta}=\frac{1}{u_{\alpha}-u_{\beta}}A_{\alpha,\beta}+\frac{1}{u_{\beta}
-u_{\alpha}}B_{\alpha,\beta}.$$ One can check that these formulas
determine a $PM$-structure on the space $\cal L$ for generic
$u_1,...,u_m,$ $ t_1,...,t_m$. The elements
$$K_{\alpha}=t_{\alpha}C_{\alpha}
+\sum_{\beta\ne\alpha}\frac{t_{\beta}}{u_{\alpha}-u_{\beta}}(A_{\alpha,\beta}B_{\beta,\alpha}-e_{\alpha})$$
satisfy the property $K_{\alpha}v=vK_{\beta}$ for all $v\in U({\cal
L})_{\alpha,\beta}$.

Note that the corresponding algebra $U(\cal L)$ has
one-dimensional representation $A_{\alpha,\beta}\to 1$,
$B_{\alpha,\beta}\to u_{\beta}/u_{\alpha}$, $C_{\alpha}\to
1/u_{\alpha},$ which gives rise to  Example 1.3.

\section{Case of semi-simple algebras $\cal A$ and $\cal B$}

\subsection{Matrix of multiplicities}

In this Subsection we suppose that $\cal L$ is a weak
$PM$-structure. We use a notation $V^l$ for a direct sum of $l$
copies of a linear space $V$ if $l\in\N$ and assume $V^0=0$. We
recall that any left $End(V)$-module has the form $V^l$ and any
right $End(V)$-module has the form $(V^{\star})^l$.

{\bf Lemma 4.1.} Let $\cal A$ be a semi-simple algebra, namely
${\cal A}=\oplus_{1\le i\le r}End(V_i),$ where $\dim V_i=m_i$.
Then $\cal L$ as $\cal A$-module is isomorphic to $\oplus_{1\le
i\le r}(V^{\star}_i)^{2m_i}$.

{\bf Proof.} It is known that any right $\cal A$-module has the form
$\oplus_{1\le i\le r}(V^{\star}_i)^{l_i}$ for some $l_1,...,l_r\ge
0$. Therefore ${\cal L}=\oplus_{1\le i\le r}{\cal L}_i$ where ${\cal
L}_i=(V^{\star}_i)^{l_i}$. Note that $\cal A\subset\cal L$ and,
moreover, $End(V_i)\subset{\cal L}_i$ for $i=1,...,r$. Besides,
$End(V_i)\bot{\cal L}_j$ for $i\ne j$. Indeed, we have
$(v,a)=(v,Id_ia)=(vId_i,a)=0$ for $v\in{\cal L}_j$ and $a\in
End(V_i),$ where $Id_i$ is the unity of the subalgebra $End(V_i)$.
Since $(\cdot,\cdot)$ is non-degenerate and $End(V_i)\bot End(V_i)$
by the property 3 of weak $PM$-structure, we have $\dim{\cal L}_i\ge
2\dim End(V_i)$. But $\sum_i\dim{\cal L}_i=\dim{\cal L}=2\dim{\cal
A}=\sum_i 2\dim End(V_i)$ and we obtain the identity $\dim{\cal
L}_i=2\dim End(V_i)$ for each $i=1,...,r,$ which is equivalent to
the statement of Lemma 4.1.

{\bf Lemma 4.2.} Let $\cal A$ and $\cal B$ be semi-simple, namely
$${\cal A}=\oplus_{1\le i\le r}End(V_i),\qquad {\cal B}=\oplus_{1\le
j\le s}End(W_j), \qquad \dim V_i=m_i, \quad\dim W_j=n_j.$$ Then
$\cal L$ as $\cal A\otimes\cal B$-module is isomorphic to
$\oplus_{1\le i\le r,1\le j\le s}(V^{\star}_i\otimes
W_j)^{a_{i,j}},$ where $a_{i,j}\ge 0$ and
\begin{equation}\label{adm}   \sum_j
a_{i,j}n_j=2m_i,\qquad \sum_i a_{i,j}m_i=2n_j.
\end{equation}

{\bf Proof.} It is known that any $\cal A\otimes\cal B$-module has
the form $\oplus_{1\le i\le r,1\le j\le s}(V^{\star}_i\otimes
W_j)^{a_{i,j}},$ where $a_{i,j}\ge 0$. Applying Lemma 4.1, we
obtain $\dim {\cal L}_i=2m_i^2,$ where ${\cal L}_i=\oplus_{1\le
j\le s}(V^{\star}_i\otimes W_j)^{a_{i,j}}$. This gives the first
equation from (\ref{adm}). The second equation can be obtained
similarly.

{\bf Definition.} The matrix $(a_{i,j})$ from Lemma 4.2 is called
matrix of multiplicities of a weak $PM$-structure $\cal L$.

{\bf Definition.} An $r\times s$ matrix $(a_{i,j})$ is called
 decomposable  if
there exist partitions $\{1,...,r\}=I\sqcup I^{\prime}$ and
$\{1,...,s\}=J\sqcup J^{\prime}$ such that $a_{i,j}=0$ for
$(i,j)\in I\times J^{\prime}\sqcup I^{\prime}\times J$.

{\bf Lemma 4.3.} If matrix of multiplicities is decomposable, then
corresponding $PM$-structure is decomposable.

{\bf Proof.} Suppose $(a_{i,j})$ is decomposable. We have $\cal
A=\cal A^{\prime}\oplus\cal A^{\prime\prime}$, $\cal B=\cal
B^{\prime}\oplus\cal B^{\prime\prime}$ and $\cal L=\cal
L^{\prime}\oplus\cal L^{\prime\prime}$ where $${\cal
A}^{\prime}=\oplus_{i\in I}End(V_i),\qquad {\cal
A}^{\prime\prime}=\oplus_{i\in I^{\prime}}End(V_i),\qquad {\cal
B}^{\prime}=\oplus_{j\in J}End(W_j),$$ $${\cal
B}^{\prime\prime}=\oplus_{j\in J^{\prime}}End(W_j), \qquad {\cal
L}^{\prime}=\oplus_{(i,j)\in I\times J}(V^{\star}_i\otimes
W_j)^{a_{i,j}},\qquad {\cal L}^{\prime\prime}=\oplus_{(i,j)\in
I^{\prime}\times J^{\prime}}(V^{\star}_i\otimes W_j)^{a_{i,j}}.$$
It is clear that this is a decomposition of $\cal L$.

{\bf Definition.} We call an $r \times s$ matrix with non-negative
integral entries $(a_{i,j})$ admissible if it is indecomposable and
(\ref{adm}) holds for some positive vectors $(m_1,...,m_r)$ and
$(n_1,...,n_s)$.

Now our aim is to classify all  admissible matrices. Note that if
$A$ is admissible, then $A^t$ is also admissible. Moreover, if $A$
is the matrix of multiplicities of a weak $PM$ structure with
semi-simple algebras $\cal A$ and $\cal B$, then $A^t$ is the matrix
of multiplicities of the opposite weak $PM$-structure.

{\bf Theorem 4.1.} There is a one-to-one correspondence between the
following two sets:

{\bf 1.} Admissible matrices up to a permutation of rows and
columns.

{\bf 2.} Simple laced affine Dynkin diagrams with a partition of
the set of vertices into two subsets (represented by black and
white circles in the pictures below) such that vertices in each
subset are pairwise non-connected.

Namely, assign to each vertex of such a Dynkin diagram a vector
space  from the set \linebreak $\{V_1,...,V_r,W_1,...,W_s\}$ in such
a way that there is a one-to-one correspondence between this set and
the set of vertices, and for any $i,j$ the spaces $V_i$, $V_j$ are
not connected by edges  as well as the spaces $W_i$, $W_j$. Then
$a_{i,j}$ is equal to the number of edges between $V_i$ and $W_j$.

{\bf Proof.} Let $(a_{i,j})$ be an admissible $r\times s$ matrix.
Consider a linear space with a basis $\{v_1,...,v_r,w_1,...,w_s\}$
and the symmetric bilinear form $(v_i,v_j)=(w_i,w_j)=2\delta_{i,j}$,
$(v_i,w_j)=-a_{i,j}$. Let $J=m_1 v_1+...+m_r v_r+n_1 w_1+...+n_s
w_s$. It is clear that the equations (\ref{adm}) can be written as
follows $(v_i,J)=(w_j,J)=0,$ which means that $J$ belongs to the
kernel of the form $(\cdot,\cdot)$. Therefore (see \cite{vinberg}),
the matrix of the form is a Cartan matrix of a simple laced affine
Dynkin diagram.

On the other hand, consider a simple laced affine Dynkin diagram
with a partition of the set of vertices into two subsets such that
vertices of the same subset are not connected. It is clear that if
such a partition exists, then it is unique up to transposition of
subsets. Let $v_1,...,v_r$ be roots corresponding to vertices of the
first subset and $w_1,...,w_s$ be roots corresponding to the second
subset. We have $(v_i,v_j)=(w_i,w_j)=2\delta_{i,j}$. Let
$a_{i,j}=-(v_i,w_j)$ and $J=m_1v_1+...+m_rv_r+n_1w_1+...+n_sw_s$ be
an imaginary root. It is clear that (\ref{adm}) holds and therefore
$(a_{i,j})$ is admissible.

Note that the transposition of the subsets corresponds to the
transposition of matrix $(a_{i,j})$.

Applying known classification of affine Dynkin diagrams \cite{burb},
we obtain the following

{\bf Theorem 4.2.} Let $A=(a_{i,j})$ be an $r\times s$ matrix of
multiplicities for a weak $PM$-structure. Then, after a possible
permutation of rows and columns and the transposition, a matrix $A$
is equal to one in the following list:

{\bf 1.} $A=(2)$. Here $r=s=1,$ $n_1=m_1=m$. The corresponding
Dynkin diagram is of the type $\tilde A_1.$

{\bf 2.} $a_{i,i}=a_{i,i+1}=1$ and $a_{i,j}=0$ for other pairs
$i,j$. Here $r=s=k\ge 2,$ the indexes are taken modulo $k$,  and
$n_i=m_i=m$. The corresponding Dynkin diagram is $\tilde A_{2
k-1}.$

{\bf 3.} $A=$\begin{math} \bordermatrix{&\cr &1&1&0&0\cr
&1&0&1&0\cr &1&0&0&1}\end{math}. Here $r=3, \, s=4$ and $n_1=3
m,\,\, n_2=n_3=n_4=m, \,$\linebreak $\,m_1=m_2=m_3=2 m$. The
Dynkin diagram is $\tilde E_{6}:$ 
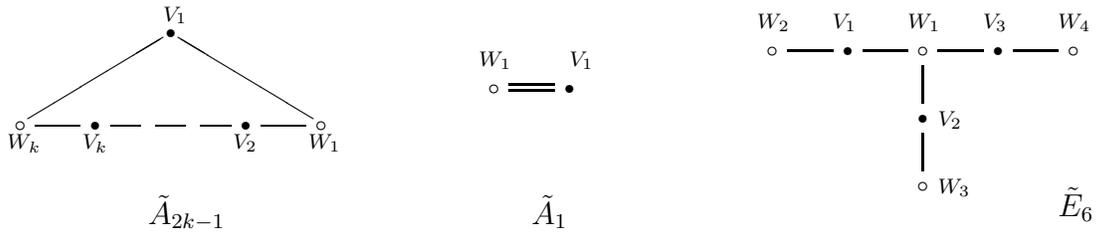
\vspace{1cm}

\begin{picture}(150,33)

\put(42,20){\circle{1.0}}

\put(32,20){\circle*{1.2}}

\put(12,20){\circle*{1.2}}

\put(2,20){\circle{1.0}}

\put(40.3,17){{\scriptsize $W_1$}}

\put(30.3,17){{\scriptsize $V_2$}}

\put(10.3,17){{\scriptsize $V_{k}$}}

\put(0.3,17){{\scriptsize $W_{k}$}}

\put(34,20){\line(1,0){6}}

\put(26,20){\line(1,0){4}}

\put(20,20){\line(1,0){4}}

\put(14,20){\line(1,0){4}}

\put(4,20){\line(1,0){6}}

\put(3,21){\line(5,3){18}}

\put(41,21){\line(-5,3){18}}

\put(22,32.3){\circle*{1.0}}

\put(21,34){{\scriptsize $V_1$}}

\put(19,7){{\bf $\tilde A_{2 k-1}$}}


\put(65,25){\circle{1.0}}

\put(75,25){\circle*{1.0}}


\put(67,25.5){\line(1,0){6}}

\put(67,24.7){\line(1,0){6}}

\put(75,28){{\scriptsize $V_1$}}

\put(63,28){{\scriptsize $W_1$}}

\put(70,7){{\bf $\tilde A_{1}$}}


\put(142,30){\circle{1.0}}

\put(140,33){{\scriptsize $W_4$}}

\put(134,30){\line(1,0){6}}

\put(124,30){\line(1,0){6}}

\put(114,30){\line(1,0){6}}

\put(104,30){\line(1,0){6}}

\put(132,30){\circle*{1.2}}

\put(130,33){{\scriptsize $V_3$}}

\put(122,30){\circle{1.0}}

\put(120,33){{\scriptsize$W_1$}}

\put(112,30){\circle*{1.2}}

\put(110,33){{\scriptsize $V_1$}}

\put(102,30){\circle{1.0}}

\put(100,33){{\scriptsize $W_2$}}

\put(122,21){\circle*{1.0}}

\put(122,12){\circle{1.0}}

\put(122,23){\line(0,1){5}}

\put(122,14){\line(0,1){5}}

\put(124,20){{\scriptsize $V_2$}}

\put(124,11){{\scriptsize $W_3$}}

\put(140,8){{\bf $\tilde E_{6}$}}

\end{picture}

{\bf 4.} $A=(1,1,1,1)$. Here $r=1,s=4$ and $n_1=n_2=n_3=n_4=m,\,\,
m_1=2 m$. The corresponding Dynkin diagram is $\tilde D_4.$

{\bf 5.} $a_{1,1}=a_{1,2}=a_{1,3}=1,\,\,$
$a_{2,3}=a_{2,4}=a_{3,4}=a_{3,5}=\cdots=a_{k-2,k-1}=a_{k-2,k}=1,\,\,$
$a_{k-1,k}=a_{k-1,k+1}=a_{k-1,k+2}=1,\,\,$
 and $a_{i,j}=0$ for other $(i,j)$.
Here we have $r=k-1,$ $s=k+2$ and $n_1=n_2=n_{k+1}=n_{k+2}=m,\,\,
n_3=\cdots=n_{k}=2 m,\,\, m_1=\cdots=m_{k}=2 m$. The corresponding
Dynkin diagram is $\tilde D_{2 k}$,\, where $k \ge 3.$

{\bf 6.} $a_{1,1}=a_{1,2}=a_{1,3}=1,\,\,$
$a_{2,3}=a_{2,4}=a_{3,4}=a_{3,5}=\cdots=a_{k-2,k-1}=a_{k-2,k}=1,\,\,$
$a_{k-1,k}=a_{k,k}=1,\,\,$
 and $a_{i,j}=0$ for other $(i,j)$.
Here we have $r=s=k\ge 3,\,\,$ $n_1=n_2=m,\,$\linebreak $
n_3=\cdots=n_{k}=2 m, \,\,m_1=\cdots=m_{k-2}=2 m,
\,\,m_{k-1}=m_k=m$. The corresponding Dynkin diagram is $\tilde
D_{2 k-1}.$ Note that if $k=3$, then
$a_{1,1}=a_{1,2}=a_{1,3}=1,\,\,$ $a_{2,3}=a_{3,3}=1$.

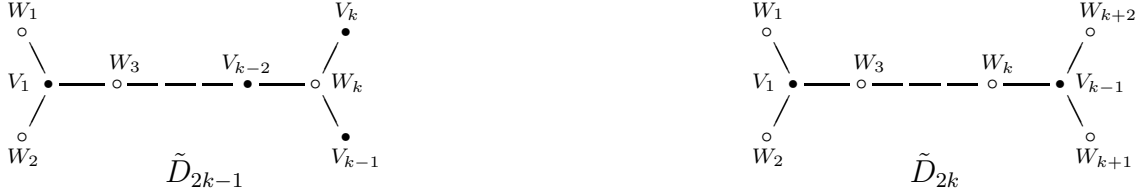
\vspace{1cm}

\begin{picture}(150,33)

\put(5.5,20){\circle*{1.0}}

\put(2,27){\circle{1.0}}

\put(0,29){{\scriptsize $W_1$}}

\put(0,9.5){{\scriptsize $W_2$}}

\put(0,19.5){{\scriptsize $V_1$}}

\put(13.5,22){{\scriptsize $W_3$}}

\put(28.5,22){{\scriptsize $V_{k-2}$}}

\put(43,19.5){{\scriptsize $W_{k}$}}

\put(43.5,29){{\scriptsize $V_{k}$}}

\put(43,9.5){{\scriptsize $V_{k-1}$}}

\put(2,13){\circle{1.0}}

\put(5,21.5){\line(-1,2){2}}

\put(5,18.5){\line(-1,-2){2}}

\put(7,20){\line(1,0){6}}

\put(16,20){\line(1,0){4}}

\put(21,20){\line(1,0){4}}

\put(26,20){\line(1,0){4}}

\put(14.6,20){\circle{1.0}}

\put(32,20){\circle*{1.0}}

\put(33.5,20){\line(1,0){6}}

\put(42,21.5){\line(1,2){2}}

\put(41,20){\circle{1.0}}

\put(42,18.5){\line(1,-2){2}}

\put(45,27){\circle*{1.0}}

\put(45,13){\circle*{1.0}}

\put(21,7){{\bf $\tilde D_{2k-1}$}}

\put(104.5,20){\circle*{1.0}}

\put(101,27){\circle{1.0}}

\put(99,29){{\scriptsize $W_1$}}

\put(99,9.5){{\scriptsize $W_2$}}

\put(99,19.5){{\scriptsize $V_1$}}

\put(112.5,22){{\scriptsize $W_3$}}

\put(129.5,22){{\scriptsize $W_{k}$}}

\put(142,19.5){{\scriptsize $V_{k-1}$}}

\put(142.5,29){{\scriptsize $W_{k+2}$}}

\put(142,9.5){{\scriptsize $W_{k+1}$}}

\put(101,13){\circle{1.0}}

\put(104,21.5){\line(-1,2){2}}

\put(104,18.5){\line(-1,-2){2}}

\put(106,20){\line(1,0){6}}

\put(115,20){\line(1,0){4}}

\put(120,20){\line(1,0){4}}

\put(125,20){\line(1,0){4}}

\put(113.6,20){\circle{1.0}}

\put(131,20){\circle{1.0}}

\put(132.5,20){\line(1,0){6}}

\put(141,21.5){\line(1,2){2}}

\put(140,20){\circle*{1.0}}

\put(141,18.5){\line(1,-2){2}}

\put(144,27){\circle{1.0}}

\put(144,13){\circle{1.0}}

\put(120,7){{\bf $\tilde D_{2k}$}}
\end{picture}

{\bf 7.} $A=$\begin{math} \bordermatrix{&\cr &1&1&0&0&0\cr
&0&1&1&1&0\cr &0&0&0&1&1}\end{math}. Here $r=3,s=5$ and
$n_1=m,\,\, n_2=3m, \,\,n_3=2m,\,$\linebreak $n_4=3m, \,\,n_5=m,
\,\,m_1=2m,\,\, m_2=4m, \,\,m_3=2m$. The Dynkin diagram is $\tilde
E_{7}.$

{\bf 8.} $A=$\begin{math} \bordermatrix{&\cr &1&0&0&0&0\cr
&1&1&1&0&0\cr &0&0&1&1&0\cr &0&0&0&1&1}\end{math}. Here $r=4,s=5$
and $n_1=4m,\,\, n_2=3m,\,\, n_3=5m,\,$\linebreak $ n_4=3m,\,\,
n_5=m,\,\, m_1=2m,\,\, m_2=6m,\,\, m_3=4m,\,\, m_4=2m$. The Dynkin
diagram is $\tilde E_{8}.$

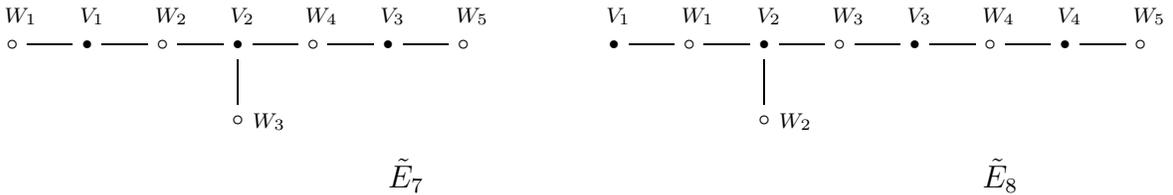
\vspace{1cm}

\begin{picture}(150,23)

\put(60,20){\circle{1.0}}

\put(50,20){\circle*{1.0}}

\put(40,20){\circle{1.0}}

\put(30,20){\circle*{1.0}}

\put(20,20){\circle{1.0}}

\put(10,20){\circle*{1.0}}

\put(0,20){\circle{1.0}}

\put(30,10){\circle{1.0}}

\put(30,12){\line(0,1){6}}

\put(52,20){\line(1,0){6}}

\put(42,20){\line(1,0){6}}

\put(32,20){\line(1,0){6}}

\put(22,20){\line(1,0){6}}

\put(12,20){\line(1,0){6}}

\put(2,20){\line(1,0){6}}

\put(49,23){{\scriptsize $V_3$}}

\put(29,23){{\scriptsize $V_2$}}

\put(9,23){{\scriptsize $V_1$}}

\put(-1,23){{\scriptsize $W_1$}}

\put(19,23){{\scriptsize $W_2$}}

\put(39,23){{\scriptsize $W_4$}}

\put(59,23){{\scriptsize $W_5$}}

\put(32,9){{\scriptsize $W_3$}}

\put(50,1){{\bf $\tilde E_{7}$}}

\put(150,20){\circle{1.0}}

\put(140,20){\circle*{1.0}}

\put(130,20){\circle{1.0}}

\put(120,20){\circle*{1.0}}

\put(110,20){\circle{1.0}}

\put(100,20){\circle*{1.0}}

\put(90,20){\circle{1.0}}

\put(80,20){\circle*{1.0}}

\put(100,10){\circle{1.0}}

\put(100,12){\line(0,1){6}}

\put(142,20){\line(1,0){6}}

\put(132,20){\line(1,0){6}}

\put(122,20){\line(1,0){6}}

\put(112,20){\line(1,0){6}}

\put(102,20){\line(1,0){6}}

\put(92,20){\line(1,0){6}}

\put(82,20){\line(1,0){6}}

\put(139,23){{\scriptsize $V_4$}}

\put(119,23){{\scriptsize $V_3$}}

\put(99,23){{\scriptsize $V_2$}}

\put(79,23){{\scriptsize $V_1$}}

\put(89,23){{\scriptsize $W_1$}}

\put(109,23){{\scriptsize $W_3$}}

\put(129,23){{\scriptsize $W_4$}}

\put(149,23){{\scriptsize $W_5$}}

\put(102,9){{\scriptsize $W_2$}}

\put(129,1){{\bf $\tilde E_{8}$}}
\end{picture}

\subsection{$PM$-structures connected with affine Dynkin diagrams}

In the previous Subsection, we have shown  that if $\cal L$ is an
indecomposable $PM$-structure with semi-simple algebras ${\cal
A}=\oplus_{1\le i\le r}End(V_i)$, ${\cal B}=\oplus_{1\le j\le
s}End(W_j)$, then there exists an affine Dynkin diagram of the type
$A$, $D$, or $E$ such that:

{\bf 1.} There is a one-to-one correspondence between the set of
vertices and the set of vector spaces $\{V_1,...,V_r,W_1,...,W_s\}$.

{\bf 2.} For any $i,j$ the spaces $V_i$, $V_j$ are not connected
by edges as well as $W_i$, $W_j$.

{\bf 3.} $\cal L$ as $\cal A\otimes\cal B$-module is isomorphic to
$\oplus_{1\le i\le r,1\le j\le s}(V^{\star}_i\otimes
W_j)^{a_{i,j}},$ where $a_{i,j}$ is equal to the number of edges
between $V_i$ and $W_j$.

{\bf 4.} The vector $(\dim V_1,...,\dim V_r,\dim W_1,...,\dim W_s)$
is an imaginary positive root of the Dynkin diagram.

To describe the corresponding $PM$-structure it remains to construct
an embedding $\cal A\to\cal L$, $\cal B\to\cal L$ and a scalar
product $(\cdot,\cdot)$ on the space $\cal L$. Note that if we fix
an element $1\in\cal L$, then we can define the embedding $\cal
A\to\cal L$, $\cal B\to\cal L$ by the formula $a\to 1a$, $b\to b1$
for $a\in\cal A$, $b\in\cal B$. After that it is not difficult to
construct a scalar product. Moreover, we may assume that $1$ is a
generic element of $\cal L$. Therefore, to study $PM$-structures
corresponding to a Dynkin diagram, one should take a generic element
in ${\cal L}=\oplus_{1\le i\le r,1\le j\le s}(V^{\star}_i\otimes
W_j)^{a_{i,j}},$ find its simplest canonical form by choosing bases
in the vector spaces $V_1,...,V_r,W_1,...,W_s$, calculate the
embedding $\cal A\to\cal L$, $\cal B\to\cal L$ and the scalar
product $(\cdot,\cdot)$ on the space $\cal L$.

For example, consider the case $\tilde A_{2 k-1}.$ We have $\dim
V_i=\dim W_i=m$ for $1\le i\le k$. Let $\{v_{i,\alpha};1\le\alpha\le
m\}$ be a basis of $V^{\star}_i$ and $\{w_{i,\alpha};1\le\alpha\le
m\}$ be a basis of $W_i$. Let
$\{e_{i,\alpha,\beta};1\le\alpha,\beta\le m\}$ be a basis of
$End(V_i)$ such that
$v_{i,\alpha}e_{i,\alpha^{\prime},\beta}=\delta_{\alpha,\alpha^{\prime}}v_{i,\beta}$
and $\{f_{i,\alpha,\beta};1\le\alpha,\beta\le m\}$ be a basis of
$End(W_i)$ such that
$f_{i,\alpha,\beta}w_{i,\beta^{\prime}}=\delta_{\beta,\beta^{\prime}}w_{i,\alpha}$.
A generic element $1\in\cal L$ in a suitable  basis in $V_i$, $W_i$
can be written in the form $1=\sum_{1\le i\le k,1\le\alpha\le
m}(v_{i,\alpha}\otimes w_{i,\alpha}+\lambda_{\alpha}
v_{i+1,\alpha}\otimes w_{i,\alpha}),$ where index $i$ is taken
modulo $k$ and $\lambda_1,...,\lambda_m\in\C$ are generic complex
numbers. The embedding $\cal A\to\cal L$, $\cal B\to\cal L$ is the
following: $e_{i,\alpha,\beta}\to
1e_{i,\alpha,\beta}=v_{i,\beta}\otimes w_{i,\alpha}+\lambda_{\alpha}
v_{i,\beta}\otimes w_{i-1,\alpha},\quad f_{i,\alpha,\beta}\to
f_{i,\alpha,\beta}1=v_{i,\beta}\otimes w_{i,\alpha}+\lambda_{\beta}
v_{i+1,\beta}\otimes w_{i,\alpha}$. It is clear that $\dim {\cal
A\cap\cal B}=m$ and a basis of this space is
$\{\sum_i(v_{i,\alpha}\otimes w_{i,\alpha}+\lambda_{\alpha}
v_{i,\alpha}\otimes w_{i-1,\alpha});1\le\alpha\le m\}$. It is also
clear that the algebra $\cal A\cap\cal B$ is isomorphic to $\C^m$.
Let us introduce a new basis in the algebras $\cal A$ and $\cal B$.
Namely, let $A^i_{\alpha,\beta}=\sum_{1\le j\le
k}\epsilon^{ij}e_{i,\alpha,\beta}$ and
$B^i_{\alpha,\beta}=\sum_{1\le j\le
k}\epsilon^{ij}f_{i,\alpha,\beta}$. Here $\epsilon=\exp(2\pi i/k)$
is a primitive root of unity of degree $k$. Simple calculations give
now the following description of the corresponding $PM$-structure in
the case $\tilde A_{2 k-1}.$

The algebra $\cal A$ has a basis $\{A^i_{\alpha,\beta};
1\le\alpha,\beta\le m, i\in \Z/k\Z\}$ such that
$A^i_{\alpha,\beta}A^j_{\beta,\gamma}=A^{i+j}_{\alpha,\gamma}.$ The
algebra $\cal B$ has a basis $\{B^i_{\alpha,\beta};
1\le\alpha,\beta\le m, i\in \Z/k\Z\}$ such that
$B^i_{\alpha,\beta}B^j_{\beta,\gamma}=B^{i+j}_{\alpha,\gamma}.$ The
intersection $\cal A\cap\cal B$ has a basis
$\{e_{\alpha}=A^0_{\alpha,\alpha}=B^0_{\alpha,\alpha};1\le\alpha\le
m\}$. A basis of the space $\cal L$ consists of the elements
$e_{\alpha}$, $A^i_{\alpha,\beta}$, $B^i_{\alpha,\beta},$ where
$i\ne 0$ if $\alpha=\beta$ and $C_{\alpha},$ where $1\le\alpha\le
m$. The scalar product has the form
$(B^i_{\alpha,\beta},A^{-i}_{\beta,\alpha})=(\epsilon^i\lambda_{\alpha}-\lambda_{\beta})/t_{\alpha}$,
$(e_{\alpha},C_{\alpha})=t_{\alpha}^{-1}$. The action of $\cal A$
and $\cal B$ on the space $\cal L$ is given by the formulas:
$$B^i_{\alpha,\beta}A^j_{\beta,\gamma}=\frac{\epsilon^{-j}\lambda_{\gamma}-\lambda_{\beta}}
{\epsilon^{-i-j}\lambda_{\gamma}-\lambda_{\alpha}}A^{i+j}_{\alpha,\gamma}+
\frac{\epsilon^i\lambda_{\alpha}-\lambda_{\beta}}{\epsilon^{i+j}\lambda_{\alpha}-\lambda_{\gamma}}B^{i+j}_{\alpha,\gamma},$$
where $i+j\ne 0$ or $\alpha\ne\gamma$  and
$$B^i_{\alpha,\beta}A^{-i}_{\beta,\alpha}=\epsilon^ie_{\alpha}+(\epsilon^i\lambda_{\alpha}-\lambda_{\beta})C_{\alpha},$$
$$C_{\alpha}A^i_{\alpha,\beta}=\frac{1}{\epsilon^{-i}\lambda_{\beta}-\lambda_{\alpha}}A^i_{\alpha,\beta}+\frac{1}
{\epsilon^i\lambda_{\alpha}-\lambda_{\beta}}B^i_{\alpha,\beta},$$
$$B^i_{\alpha,\beta}C_{\beta}=\frac{1}{\epsilon^{-i}\lambda_{\beta}-\lambda_{\alpha}}A^i_{\alpha,\beta}+\frac{1}
{\epsilon^i\lambda_{\alpha}-\lambda_{\beta}}B^i_{\alpha,\beta}.$$
Here $\epsilon=\exp(2\pi i/k)$ and
$\lambda_1,...,\lambda_m,t_1,...,t_m\in\C$ such that
$(\lambda_{\alpha})^k\ne(\lambda_{\beta})^k$ for $\alpha\ne\beta$
and $t_{\alpha}\ne 0$. The elements
$$K_{\alpha}=t_{\alpha}C_{\alpha}+\sum_{(i,\beta)\ne(0,\alpha)}\frac{t_{\beta}}{\epsilon^i\lambda_{\beta}-
\lambda_{\alpha}}(A^{-i}_{\alpha,\beta}B^i_{\beta,\alpha}-\epsilon^ie_{\alpha})$$
satisfy the property $K_{\alpha}v=vK_{\beta}$ for all $v\in U({\cal
L})_{\alpha,\beta}$.

The corresponding operator $R$ has the following components:
$$R_{\beta,\alpha}(x_{\alpha})=\sum_{i\in\Z/k\Z}
\frac{t_{\alpha}}{\epsilon^i\lambda_{\alpha}-\lambda_{\beta}}a_{\beta,\alpha}^{-i}x_{\alpha}b^i_{\alpha,\beta}$$
for $\alpha\ne\beta$ and
$$R_{\alpha,\alpha}(x_{\alpha})=t_{\alpha}c_{\alpha}x_{\alpha}+
\sum_{(i,\beta)\ne(0,\alpha)}\frac{t_{\beta}}{\epsilon^i\lambda_{\beta}-
\lambda_{\alpha}}(a^{-i}_{\alpha,\beta}x_{\alpha}b^i_{\beta,\alpha}-\epsilon^ix_{\alpha}).$$
Here $A^i_{\alpha,\beta}\to a^i_{\alpha,\beta}$,
$B^i_{\alpha,\beta}\to b^i_{\alpha,\beta}$ and $C_{\alpha}\to
c_{\alpha}$ is a representation.

Let $a$, $t$ be linear operators in some vector space. Assume that
$a^k=1$, $at=\epsilon ta$ and the operators $t-\lambda_{\alpha}$ are
invertible for $1\le\alpha\le m$. It is easy to check that the
formulas
$$A_{\alpha,\beta}^i\to a^i,\qquad B_{\alpha,\beta}^i\to
\frac{\epsilon^i t-\lambda_{\beta}}{t-\lambda_{\alpha}}a^i,\qquad
C_{\alpha}\to\frac{1}{t-\lambda_{\alpha}}$$ define a representation
of the algebra $U(\cal L)$. Note that we do not assume that $t^k=1$.
We have only $at^k=t^ka$ which easily follows from the commutation
relation between $a$ and $t$.

{\bf Remark 1.} If $m=1$, then this is the Example 2.1. If $k=1$,
then this is the Example 3.1.

{\bf Remark 2.} Since operator $R$ depends linearly on
$t_1,...,t_m$, we obtain $m+1$ pairwise compatible multiplications.
One can check that these multiplications can be obtained using
Theorem 1.1 from the case $m=1$. We conjecture that the similar
result holds for other Dynkin diagrams.

The cases corresponding to affine Dynkin diagrams of type $D$ and
$E$ are treated similarly, but resulting formulas are more
complicated. Note that classification of generic elements $1\in\cal
L$ up to choice of bases in the vector spaces
$V_1,...,V_r,W_1,...,W_s$ is equivalent to classification of
representations of a quiver corresponding to our affine Dynkin
diagram with the same vector spaces. Therefore, we can apply known
results about these representations (see \cite{quivers1, quivers2,
quivers3}). Since the dimension of a representation is equal to $m
I,$ where $I$ is the minimal positive imaginary root and our
representation is generic, then it is isomorphic to a direct sum of
$m$ irreducible representations of dimension $I$. Therefore,
$1=e_1+...+e_m,$ where $e_1,...,e_m$ correspond to these
representations. Taking  the explicit form of these representations
for affine Dynkin diagrams of the type $D$ and $E$ from
\cite{quivers2, quivers3} and applying the scheme described above,
one can obtain explicit formulas for the corresponding
$PM$-structures similarly to the case of the diagrams of the type
$A$.

\vspace{8mm} \centerline{\Large\bf Conclusion}
\medskip

In this paper we have studied associative multiplications in a
semi-simple associative algebra over $\C$ compatible with the usual
one. It turned out that these multiplications are in one-to-one
correspondence with representations of  $M$-structures in the matrix
case and $PM$-structures in the case of direct sum of several matrix
algebras. These structures are differ from the Hopf algebras but in
some features remind them. Namely, a $PM$-structure also contains
two (associative) algebras $\cal A$ and $\cal B,$ which are dual in
some sense and satisfy certain compatibility conditions between
them. Natural problem arises: to classify $PM$-structures for
semi-simple algebras $\cal A$ and $\cal B$ (this is done in Section
{\bf 4}) or, which is more difficult, to describe $PM$-structures if
only one of these algebras is semi-simple (the case of commutative
semi-simple $\cal A$ is treated in Subsection {\bf 2.3}).

Another interesting question is to investigate integrable systems
corresponding to given representations of $PM$-structures. The
problem here is to formulate properties of the integrable system in
terms of algebraic properties of $PM$-algebra. It would be also
interesting to study corresponding quantum integrable systems.

Note that our $M$ and $PM$-structures are the particular cases of
the following general situation. We have a linear space $\cal L$
with two subspaces $\cal A$ and $\cal B$ and a non-degenerate scalar
product. The spaces $\cal A$ and $\cal B$ are associative algebras
with common subalgebra $\cal S=\cal A\cap\cal B$. We assume that
$\dim \cal S=\dim\cal A\cap\cal B=\dim\cal L/(\cal A+\cal B)$ and
our scalar product restricted on $\cal A$ and on $\cal B$ is zero.
We have also a left action of $\cal A$ and a right action of $\cal
B$ on $\cal L$ which commute with each other and invariant with
respect to the scalar product (that is $(b_1b_2,v)=(b_1,b_2v)$,
$(v,a_1a_2)=(va_1,a_2)$ for any $a_1,a_2\in\cal A$, $b_1,b_2\in\cal
B$ and $v\in\cal L$). Finally, we assume that $\cal A\subset\cal L$
is a submodule with respect to the action of $\cal A,$ where $\cal
A$ acts by right multiplication and similar property is valid for
$\cal B$. Now, if ${\cal S}=0$, then we have the toy example from
the Introduction, if $\cal S=\C$, then we have a weak $M$-structure
and if $\cal S$ is a direct sum of $m$ copies of $\C$, then we have
a weak $PM$-structure of size $m$. It would be interesting to study
and find possible applications of these structures for different
$\cal S$.

\vskip.3cm \noindent {\bf Acknowledgments.} The authors are
grateful to I.Z. Golubchik for useful discussions. The research
was partially supported by: RFBR grant 05-01-00189, NSh grants
1716.2003.1 and 2044.2003.2.

\end{document}